\documentclass[a4paper,12pt]{article}
\begin{document}

\title{Left invariant measures on locally compact fan loops.}
\author{S.V. Ludkowski}
\date{10 December 2018}
\maketitle
\begin{abstract}
In this article left invariant measures and functionals on locally
compact nonassociative fan loops are investigated. For this purpose
necessary properties of topological fan loops, estimates and
approximations of functions on them are studied. An existence of
nontrivial left invariant measures on locally compact fan loops is
proved. Abundant families of fan loops are provided with the help of
different types of their products. \footnote{key words and phrases:
measure; left invariant;
loop; locally compact  \\
Mathematics Subject Classification 2010: 43A05; 28C10; 20N05; 22A30}

\end{abstract}
\par Address: Dep. Appl. Mathematics, Moscow State Techn. Univ. MIREA,
\par av. Vernadsky 78, Moscow 11944, Russia; e-mail: sludkowski@mail.ru

\section{Introduction.}
Left invariant measures or Haar measures on locally compact groups
play very important role in measure theory, harmonic analysis,
representation theory, geometry, mathematical physics, etc. (see,
for example, \cite{boloktodb,fell,hew} and references therein). On
the other hand, in nonassociative algebra, in noncommutative
geometry, field theory, topological algebra there frequently appear
binary systems which are nonassociative generalizations of groups
and related with loops, quasi-groups, Moufang loops, IP-loops, etc.
(see \cite{bruckb,kakkar,kiechlb,pickert,razm,smithb,vojtech} and
references therein). An arbitrary IP-loop $Y$ is a loop with a
restriction: for each $x\in Y$ there exist elements $x_1$ and $x_2$
in $Y$ such that for each $y$ in $Y$ the identities are satisfied
$x_1(xy)=y$ and $(yx)x_2=y$, where $x_1$ and $x_2$ are also denoted
by $\mbox{}^{-1}x$ and $x^{-1}$ and called left and right inverses
of $x$ respectively. It was investigated and proved in the 20-th
century that a nontrivial geometry exists if and only if there
exists a corresponding loop.
\par Very important role in mathematics and quantum field theory play
octonions and generalized Cayley-Dickson algebras
\cite{albmajja99,allcja98,baez,dickson}. A multiplicative law of
their canonical bases is nonassociative and leads to a more general
notion of a metagroup instead of a group \cite{ludlmla18}. They are
used not only in algebra and geometry, but also in noncommutative
analysis and PDEs, particle physics, mathematical physics (see
\cite{baez,castdoyfior,dickson,frenludkfejms18}-\cite{guetze,kansol}-\cite{ludkcvee13}
and references therein). The preposition "meta" is used to emphasize
that such an algebraic object has properties milder than a group. By
their axiomatic metagroups are loops with weak relations. They were
used in \cite{ludlmla18} for investigations of automorphisms and
derivations of nonassociative algebras. \par In this article more
general binary systems such as fan loops are studied (see Definition
2.1). They also are more general than IP-loops, because in fan loops
$G$ left and right inverses $\mbox{}^{-1}x$ and $x^{-1}$ of nonunit
elements $x$ in $G$ may not exist.
\par This article is devoted to left invariant measures (see
Definition 3.18) on locally compact fan loops. Necessary preliminary
results about fan loops are given in Section 2. Specific algebraic
and topological features of fan loops are studied in Lemmas 2.2-2.6,
2.12 and Propositions 2.9, 2.11. A quotient of a fan loop by its fan
is investigated in Theorem 2.8. A uniform continuity of maps on
topological fan loops is studied in Theorem 2.14 and Corollary 2.15.
\par Left invariant functionals and measures are investigated in
Section 3. These properties are more complicated than for groups and
IP-loops, because of the nonassociativity of fan loops and absence
of left and right inverses in general. Main results are theorems
3.15, 3.16, 3.19, 3.20. For their proofs estimates of nonnegative
functions with compact supports in fan loops are investigated in
Lemmas 3.2, 3.4, 3.6. Functionals on a space of nonnegative
functions with compact supports in a fan loop are studied in Lemmas
3.7, 3.8, 3.10, 3.13 and Theorem 3.9. In Theorem 3.11 approximations
of nonnegative functions with compact supports in the fan loop are
described. \par In an appendix abundant families of fan loops are
provided with the help of a direct product and smashing products
(see Remark 4.3 and Definition 4.5). For this purpose Theorems 4.1
and 4.4 are proved.
\par All main results of this paper are obtained for the first time.
They can be used in harmonic analysis on nonassociative algebras and
metagroups and loops, representation theory, geometry, mathematical
physics, quantum field theory, particle physics, PDEs, etc.

\section{Fan loops.}
\par To avoid misunderstandings we give necessary definitions.
For short it will be written fan loop instead of nonassociative fan
loop.
\par {\bf 2.1. Definition.}  Let $G$ be a set with a multiplication
(that is a single-valued binary operation) $G^2\ni (a,b)\mapsto ab
\in G$ defined on $G$ satisfying the conditions:
\par $(2.1.1)$ for each $a$ and $b$ in $G$ there is a unique $x\in
G$ with $ax=b$ and
\par $(2.1.2)$ a unique $y\in G$ exists satisfying $ya=b$, which are
denoted by $x=a\setminus b=Div_l(a,b)$ and $y=b/a=Div_r(a,b)$
correspondingly,
\par $(2.1.3)$ there exists a neutral (i.e. unit) element $e_G=e\in G$: $~eg=ge=g$
for each $g\in G$.
\par We consider subsets in $G$:
\par $(2.1.4)$ $Com (G) := \{ a\in G: ~ \forall b\in G, ~ ab=ba \} $;
\par $(2.1.5)$ $N_l(G) := \{a\in G: ~ \forall b\in G, ~ \forall c\in G, ~ (ab)c=a(bc) \}
$;
\par $(2.1.6)$ $N_m(G) := \{a\in G: ~ \forall b\in G, ~ \forall c\in G, ~ (ba)c=b(ac)
\} $;
\par $(2.1.7)$ $N_r(G) := \{a\in G: ~ \forall b\in G, ~ \forall c\in G, ~ (bc)a=b(ca)
\} $;
\par $(2.1.8)$ $N(G) := N_l(G)\cap N_m(G)\cap N_r(G)$;
\par $Z(G) := Com (G)\cap N(G)$.
\par Then $N(G)$ is called a nucleus of $G$
and $Z(G)$ is called the center of $G$.
\par We call $G$ a fan loop if a set $G$ possesses a multiplication
and satisfies conditions $(2.1.1)$-$(2.1.3)$ and
\par $(2.1.9)$ $(ab)c=t(a,b,c)a(bc)$ and $(ab)c=a(bc)p(a,b,c)$\\ for each
$a$, $b$ and $c$ in $G$, where \par $t(a,b,c)=t_G(a,b,c)\in N(G)$
and $p(a,b,c)=p_G(a,b,c)\in N(G)$.
\par Then $G$ will be called a central fan loop if in addition to $(2.1.9)$ it
satisfies the condition:
\par $(2.1.10)$ $ab={\sf t}_2(a,b)ba$ \\ for each $a$ and $b$ in $G$, where
${\sf t}_2(a,b)\in Z(G)$.
\par Let $\tau $ be a topology on $G$ such that
the multiplication $G\times G\ni (a,b)\mapsto ab\in G$ and the
mappings $Div_l(a,b)$ and $Div_r(a,b)$ are jointly continuous
relative to $\tau $, then $(G, \tau )$ will be called a topological
fan loop. Henceforth it will be assumed that $\tau $ is the $T_1\cap
T_{3.5}$ topology, if something other will not be specified.
\par A minimal closed subgroup $N_0(G)$ in the topological fan loop $G$ containing $t(a,b,c)$ and
$p(a,b,c)$ for each $a$, $b$ and $c$ in $G$ will be called a fan of
$G$.

\par Elements of the fan loop $G$ will be denoted by small letters,
subsets of $G$ will be denoted by capital letters. If $A$ and $B$
are subsets in $G$, then $A-B$ means the difference of them $A-B=\{
a\in A: ~a \notin B \} $. Henceforward, maps and functions on fan
loops are supposed to be single-valued if something other will not
be specified.

\par {\bf 2.2. Lemma.} {\it If $G$ is a fan loop, then
for each $a$, $b$ and $c$ in $G$ the following identities are
fulfilled:
\par $(2.2.1)$ $ ~ b\setminus e=t(e/b,b,b\setminus e)(e/b)$;
\par $(2.2.1')$ $b\setminus e=(e/b)p(e/b,b,b\setminus e)$;
\par $(2.2.2)$ $(a\setminus e)b=t(e/a,a,a\setminus e)[t(e/a,a,a\setminus b)]^{-1}(a\setminus b)$;
\par $(a\setminus b)= (a\setminus e)bp(a,a\setminus e,b)$;
\par $(2.2.2')$ $(bc)\setminus a=(c\setminus (b\setminus
a))[p(b,c,(bc)\setminus a)]^{-1}$;
\par $(2.2.2'')$ $(a\setminus b)c=(a\setminus (bc))[p(a,a\setminus b,c)]^{-1}$;
\par $(2.2.2''')$ $(ab)\setminus e = (b\setminus e)(a\setminus
e)[t(a,b,b\setminus e)]^{-1}t(ab,b\setminus e,a\setminus e)$;
\par $(2.2.3)$ $b(e/a)=(b/a)p(b/a,a,a\setminus e)[p(e/a,a,a\setminus e)]^{-1} $;
\par $(b/a)= [t(b,e/a,a)]^{-1}b(e/a)$;
\par $(2.2.3')$ $a/(bc)=t(a/(bc),b,c)((a/c)/b)$;
\par $(2.2.3'')$ $c(b/a)=t(c,b/a,a)(cb)/a$;
\par $(2.2.3''')$ $e/(ab)=[p(e/b,e/a,ab)]^{-1}p(e/a,a,b)(e/b)(e/a)$.}
\par {\bf Proof.} Note that $N(G)$ is a subgroup in $G$ due to Conditions
$(2.1.5)$-$(2.1.8)$ (see also \cite{bruckb}). Then Conditions
$(2.1.1)$-$(2.1.3)$ imply that
\par $(2.2.4)$ $b(b\setminus a)=a$, $~b\setminus (ba)=a$;
\par $(2.2.5)$ $(a/b)b=a$, $~(ab)/b=a$ \\
for each $a$ and $b$ in any loop $G$ (see also
\cite{bruckb,smithb}). Using Condition $(2.1.9)$ and Identities
$(2.2.4)$ and $(2.2.5)$ we deduce that
\par $e/b=(e/b)(b(b\setminus e))  = [t(e/b,b,b\setminus e)]^{-1}(b\setminus e)
$\\ which leads to $(2.2.1)$.
\par Let $c=a\setminus b$, then from Identities $(2.2.1)$ and $(2.2.4)$ it
follows that \par $(a\setminus e)b=t(e/a,a,a\setminus
e)(e/a)(ac)$\par $=t(e/a,a,a\setminus e)[t(e/a,a,a\setminus b)]^{-1}
((e/a)a)(a\setminus b)$\\
which taking into account $(2.2.5)$ provides $(2.2.2)$.
\par On the other hand, $b\setminus e=((e/b)b)(b\setminus e) =
(e/b)(b(b\setminus e))p(e/b,b,b\setminus e)$ that gives $(2.2.1')$.
\par Let now $d=b/a$, then Identities $(2.2.1')$ and $(2.2.5)$ imply
that  \par $b(e/a)=(da)(a\setminus e)[p(e/a,a,a\setminus
e)]^{-1}$\par $=(b/a)p(b/a,a,a\setminus e)[p(e/a,a,a\setminus e)]^{-1} $ \\
which demonstrates $(2.2.3)$.
\par Next we infer from $(2.1.9)$ and $(2.2.4)$ that
\par $b(c((bc)\setminus a))=(bc)((bc)\setminus a)[p(b,c,(bc)\setminus
a)]^{-1}= a[p(b,c,(bc)\setminus a)]^{-1}$, hence $c((bc)\setminus
a)=(b\setminus a)[p(b,c,(bc)\setminus a)]^{-1}$ that implies
$(2.2.2')$.
\par Symmetrically it is deduced that
$(a/(bc))b)c=t(a/(bc),b,c)a$, consequently,
$(a/(bc))b=t(a/(bc),b,c)(a/c)$. From the latter identity it follows
$(2.2.3')$. \par Evidently, formulas \par $a((a\setminus
b)c)=(a(a\setminus b))c[p(a,a\setminus b,c)]^{-1}=bc[p(a,a\setminus
b,c)]^{-1}$ and \par $(c(b/a))a=t(c,b/a,a)cb$ \\ imply $(2.2.2'')$
and $(2.2.3'')$ correspondingly.
\par From $(2.1.9)$ we infer that \par $(ab)((b\setminus
e)(a\setminus e))=[t(ab,b\setminus e,a\setminus
e)]^{-1}t(a,b,b\setminus e)$, since by $(2.2.4)$
\par $(a(b(b\setminus e)))(a\setminus e)=e$. \\ This together with
$(2.1.1)$ and $(2.1.2)$ implies $(2.2.2''')$.
\par Analogously form $(2.1.9)$ we deduce that
\par $((e/b)(e/a))(ab)=[p(e/a,a,b)]^{-1}p(e/b,e/a,ab)$, since by
$(2.2.5)$
\par $(e/b)(((e/a)a)b)=e$.
\\ Finally applying $(2.1.1)$ and $(2.1.2)$ we get Identity $(2.2.3''')$.

\par {\bf 2.3. Lemma.} {\it Assume that $G$ is a fan loop.
Then for every $a$, $a_1$, $a_2$, $a_3$ in $G$ and $z_1$, $z_2$,
$z_3$ in $Z(G)$, $b\in N(G)$:
\par $(2.3.1)$ $t(z_1a_1,z_2a_2,z_3a_3)=t(a_1,a_2,a_3)$;
\par $(2.3.1')$ $p(z_1a_1,z_2a_2,z_3a_3)=p(a_1,a_2,a_3)$;
\par $(2.3.2)$ $t(a,a\setminus e,a)a=ap(a,a\setminus e,a)$;
\par $(2.3.2')$ $t(a,e/a,a)a=ap(a,e/a,a)$;
\par $(2.3.2'')$ $p(a,a\setminus e,a)t(e/a,a,a\setminus e)=e$;
\par $(2.3.3)$ $t(a_1,a_2,a_3b)=t(a_1,a_2,a_3)$;
\par $(2.3.3')$ $p(ba_1,a_2,a_3)=p(a_1,a_2,a_3)$;
\par $(2.3.4)$ $t(ba_1,a_2,a_3)=bt(a_1,a_2,a_3)b^{-1}$;
\par $(2.3.4')$ $p(a_1,a_2,a_3b)=b^{-1}p(a_1,a_2,a_3)b$.}
\par {\bf Proof.} Since $(a_1a_2)a_3=t(a_1,a_2,a_3)a_1(a_2a_3)$ and
$t(a_1,a_2,a_3)\in N(G)$ for every $a_1$, $a_2$, $a_3$ in $G$, then
\par $(2.3.5)$ $t(a_1,a_2,a_3)=((a_1a_2)a_3)/(a_1(a_2a_3))$. \par Therefore, for
every $a_1$, $a_2$, $a_3$ in $G$ and $z_1$, $z_2$, $z_3$ in $Z(G)$
we infer that \par $t(z_1a_1,z_2a_2,z_3a_3)=
(((z_1a_1)(z_2a_2))(z_3a_3))/((z_1a_1)((z_2a_2)(z_3a_3)))$\par $=
((z_1z_2z_3)((a_1a_2)a_3))/((z_1z_2z_3)(a_1(a_2a_3)))=
((a_1a_2)a_3)/(a_1(a_2a_3))$, since \par $(2.3.6)$
$b/(qa)=q^{-1}b/a$ and $b/q=q\setminus b=bq^{-1}$ \\ for each $q\in
Z(G)$, $ ~ a$ and $b$ in $G$, because $Z(G)$ is the commutative
group satisfying Conditions $(2.1.4)$ and $(2.1.8)$. Thus
$t(z_1a_1,z_2a_2,z_3a_3)=t(a_1,a_2,a_3)$.
\par Symmetrically we get \par $(2.3.7)$ $p(a_1,a_2,a_3)=(a_1(a_2a_3))\setminus ((a_1a_2)a_3)$ and
\par $p(z_1a_1,z_2a_2,z_3a_3)=
((z_1a_1)((z_2a_2)(z_3a_3)))\setminus
(((z_1a_1)(z_2a_2))(z_3a_3))$\par
$=((z_1z_2z_3)(a_1(a_2a_3)))\setminus ((z_1z_2z_3)((a_1a_2)a_3))=
(a_1(a_2a_3))\setminus ((a_1a_2)a_3)$ \\ that provides $(2.3.1')$.
\par From Formulas $(2.3.5)$ and $(2.2.1)$ it follows that \par $t(a,a\setminus
e,a)=((a(a\setminus e))a)/(a((a\setminus
e)a))=a/[at(e/a,a,a\setminus e)]$, consequently, \par
\par $(2.3.8)$ $t(a,a\setminus e,a)at(e/a,a,a\setminus e)=a$. \\ Then
from Formulas $(2.3.7)$, $(2.2.4)$ and Condition $(2.1.9)$ we deduce
that \par $p(a,a\setminus e,a)=(a((a\setminus e)a))\setminus
((a(a\setminus e))a)= \{ [t(a,a\setminus e,a)]^{-1}a \} \setminus
a$,
\\ which implies $(2.3.2)$. Identities $(2.3.2)$ and $(2.3.8)$ lead to
$(2.3.2'')$. Next using $(2.3.7)$ and $(2.1.9)$ we deduce that \par
$p(a,e/a,a)=[a((e/a)a)]\setminus [(a(e/a))a]=a\setminus [t(a,e/a,a)a
]$ \\ that implies $(2.3.2')$. From $(2.1.9)$ we get that
\par $((a_1a_2)a_3)b=(a_1a_2)(a_3b)=(t(a_1,a_2,a_3b)a_1(a_2a_3))b$, \\
from which and $(2.2.5)$ and $(2.3.5)$ Identity $(2.3.3)$ follows,
because $b\in N(G)$. Then
\par $b((a_1a_2)a_3)=((ba_1)a_2)a_3=b(a_1(a_2a_3)p(ba_1,a_2,a_3))$ \\
and $(2.2.4)$ and $(2.3.7)$ imply Identity $(2.3.3')$. Symmetrically
we deduce
\par $b((a_1a_2)a_3)=t(ba_1,a_2,a_3))b(a_1(a_2a_3))$ and
\par $((a_1a_2)a_3)b=(a_1(a_2a_3))bp(a_1,a_2,a_3b)$ \\ that together with
$(2.3.5)$ and $(2.3.7)$ imply Identities $(2.3.4)$ and $(2.3.4')$.

\par {\bf 2.4. Lemma.} {\it If $(G, \tau )$ is a topological loop, then
the functions $t(a_1,a_2,a_3)$ and $p(a_1,a_2,a_3)$ are jointly
continuous in $a_1$, $a_2$, $a_3$ in $G$.}
\par {\bf Proof.} This follows immediately from Formulas $(2.3.5)$, $(2.3.7)$ and
Definition 2.1.

\par {\bf 2.5. Lemma.} {\it Assume that $(G,\tau )$ is a topological
loop and $U$ is an open subsets in $G$, then for each $b\in G$ sets
$Ub$ and $bU$ are open in $G$.}
\par {\bf Proof.} Take any $c\in Ub$ and consider the equation
\par $(2.5.1)$ $xb=c$. \par Then from $(2.1.2)$ it follows that \par $(2.5.2)$ $x=c/b$.
\par Thus $x=\psi _b(c)$, where $\psi _b(c)=c/b$
is a continuous bijective function in the variable $c$ due to
Identity $(2.2.3)$ and Lemma
2.4. On the other hand, the right shift mapping \par $(2.5.3)$ $R_bu:=ub$ \\
from $G$ into $G$ is continuous and bijective in $u$ (see Definition
1). Moreover, $\psi _b(R_bu)=u$ and $R_b(\psi _b(c))=c$ for each
fixed $b\in G$ and all $u\in G$ and $c\in G$ by Identities
$(2.2.5)$. Thus $R_b$ and $\psi _b$ are open mappings, consequently,
$Ub$ is open in $G$.
\par Similarly for the equation
\par $(2.5.4)$ $by=c$ the unique solution is
\par $(2.5.5)$ $y=b\setminus c$
by Condition $(2.1.1)$. \par Therefore, $y=\theta _b(c)$, where
$\theta _b(c)=b\setminus c$ is a continuous bijective function in
$c$ according to Lemma 2.4 and Formula $(2.2.2)$. Next we consider
the left shift mapping
\par $(2.5.6)$ $L_bu=bu$ \\ for each fixed $b\in G$ and any $u\in G$.
This mapping $L_b$ is continuous, since the multiplication on $G$ is
continuous. Then $L_b(\theta _b(c))=c$ and $\theta _b(L_bu)=u$ for
every fixed $b\in G$ and all $u\in G$ and $c\in G$ by Identities
$(2.2.4)$. Therefore $\theta _b$ and $L_b$ are open mappings. Thus
the subset $bU$ is open in $G$.

\par {\bf 2.6. Lemma.} {\it Let $(G,\tau )$ be a topological
loop. \par $(i)$. Let also $U$ and $V$ be subsets in $G$ such that
either $U$ or $V$ is open, then $UV$ is open in $G$.
\par $(ii)$. If $A$ and $B$ are compact subsets in $G$, then $AB$
is compact.
\par $(iii)$. For each open neighborhood $U$ of $e$ in $G$ there
exists an open neighborhood $V$ of $e$ such that
\par $(2.6.1)$  $\check{V}\subseteq U$, where
\par $(2.6.2)$ $\check{V} =V\cup Inv_l(V)\cup Inv_r(V)$, \\ where $Inv_l(a)=Div_l(a,e)$,
$~Inv_r(a)=Div_r(a,e)$ for each $a\in G$,
\par $(2.6.3)$ $DQ=\{ x=ab: ~ a\in D, ~ b \in Q \} $,
\par $(2.6.4)$ $Inv_l(D) = \{ x=a\setminus e: ~ a \in D \} $,
\par $(2.6.5)$ $Inv_r(D) = \{ x=e/a: ~ a \in D \} $ for any subsets $D$ and $Q$ in $G$.}
\par {\bf Proof.} $(i)$. In view of Lemma 2.5 $Ub$ and $aV$ are open in $G$ for each $a\in U$
and $b\in V$, consequently, $UV=\{ x: ~ x=uv, ~ u\in U, ~ v\in V
\}=\bigcup_{b\in V}Ub= \bigcup_{a\in U}aV$ is open in $G$.
\par $(ii)$. A subset $AB = \{ c: ~ c=ab, ~ a\in A, ~ b\in B \} $ is a
continuous image of a compact subset $A\times B$ in $G\times G$,
where $G\times G$ is supplied with the product (i.e. Tychonoff)
topology, consequently, $AB$ is a compact subset in $G$ (see Theorem
3.1.10 and the Tychonoff Theorem 3.2.4 in \cite{eng}).
\par $(iii)$. The mappings $Inv_l$ and $Inv_r$ are homeomorphisms of $G$ onto
itself as the topological space, since they are bijective,
continuous and
\par $(2.6.6)$ $Inv_l(Inv_r(b))=b$ and $Inv_r(Inv_l(b))=b$ \\ for
each $b$ in $G$ by $(2.2.4)$, $(2.2.5)$. Therefore for each open
neighborhood $U$ of $e$ there exists an open neighborhood of $e$ of
the form \par $(2.6.7)$ $V:=\hat{U}$, where $\hat{U}:=U\cap
Inv_l(U)\cap Inv_r(U)$.
\par From $(2.6.6)$ we infer that $Inv_r(Inv_l(U))=U$ and
$Inv_l(Inv_r(U))=U$, hence $Inv_l(V)\subseteq U\cap Inv_l(U)\cap
Inv_l(Inv_l(U))\subseteq U\cap Inv_l(U)$ and $Inv_r(V)\subset U\cap
Inv_r(U)$, consequently, $V\cup Inv_l(V)\cup Inv_r(V)\subseteq U$.

\par {\bf 2.7. Definition.} A subloop $H$ of a loop $G$ is called
normal if it satisfies
\par $(2.7.1)$ $xH=Hx$ and
\par $(2.7.2)$ $(xy)H=x(yH)$ and $(xH)y=x(Hy)$ and $H(xy)=(Hx)y$ \\ for each $x$ and $y$
in $G$.
\par A family of cosets $\{ bH: ~ b\in G \} $ will be denoted by $G/\cdot
/N_0$.

\par {\bf 2.8. Theorem.} {\it If $G$ is a $T_1$ topological fan loop,
then its fan $N_0$ is a normal subgroup and its quotient $G/\cdot
/N_0$ is a $T_1\cap T_{3.5}$ topological group.}
\par {\bf Proof.} Let $\tau $ be a $T_1$ topology on $G$ relative to
which $G$ is a topological loop. Then each point $x$ in $G$ is
closed, since $G$ is the $T_1$ topological space (see Section 1.5 in
\cite{eng}). From the joint continuity of the multiplication and the
mappings $Div_l$ and $Div_r$ it follows that the nucleus $N=N(G)$ is
closed in $G$. Therefore the subgroup $N_0$ is the closure of a
subgroup $N_{0,0}(G)$ in $N$ generated by elements $t(a,b,c)$ and
$p(a,b,c)$ for all $a$, $b$ and $c$ in $G$ (see Definition 2.1).
According to $(2.1.5)$-$(2.1.8)$ one gets that $N$ and hence $N_0$
are subgroups in $G$ satisfying Conditions $(2.7.2)$, because
$N_0\subseteq N$ (see also \cite{bruckb,smithb}).
\par Let $a$ and $b$ belong to $N$ and $x\in G$.
Then $x(x\setminus (ab))=ab$ and \\ $x((x\setminus
a)b)=(x(x\setminus a))b=ab$, consequently, \par $(2.8.1)$
$x\setminus (ab)=(x\setminus a)b$ for each $a$ and $b$ in $N(G)$,
$x\in G$. \\ Similarly it is deduced \par $(2.8.2)$ $(ab)/x=a(b/x)$
for each $a$ and $b$ in $N(G)$, $x\in G$. \par Therefore from
$(2.1.9)$ and $(2.2.4)$ and $(2.8.1)$ it follows that
\par $((x\setminus a)x)((x\setminus b)x)=(x\setminus
a)(x((x\setminus b)x))p(x\setminus a,x,(x\setminus b)x)$\par
$=(x\setminus (ab))x[p(x,x\setminus b,x)]^{-1}p(x\setminus
a,x,(x\setminus b)x)$,
\\ since $(x\setminus a)(bx)=((x\setminus a)b)x=(x\setminus (ab))x$.
Thus \par $(2.8.3)$ $(x\setminus (ab))x= ((x\setminus
a)x)((x\setminus b)x)[p(x\setminus a,x,(x\setminus b)x)]^{-1}
p(x,x\setminus b,x)$ for each $a$ and $b$ in $N(G)$, $x\in G$.
\par From Identities $(2.2.2')$ and $(2.2.2'')$ it follows that
\par $(2.8.4)$ $x\setminus ((u\setminus v)y)=((ux)\setminus (vy))p(u,x,(ux)\setminus (vy)) [p(u,u\setminus v,x)]^{-1}$
\\  for each $u$, $v$, $x$ and $y$ in $G$, since
\par $x\setminus ((u\setminus v)y)=x\setminus (u\setminus (vy))[p(u,u\setminus
v,y)]^{-1}$. \par In particular for $u=a(bc)$ and $v=(ab)c$ with any
$a$, $b$ and $c$ in $G$ we infer using $(2.1.9)$ that
$ux=(a(b(cx)))p(b,c,x)p(a,bc,x)$ and $vx=(ab)(cx)p(ab,c,x)$, hence
from $(2.8.4)$ and $(2.3.7)$ it follows that \par $(2.8.5)$
$x\setminus
(p(a,b,c)x)=[p(b,c,x)p(a,bc,x)]^{-1}p(a,b,cx)p(u,x,(ux)\setminus
(vx))$, since \par $x\setminus
(p(a,b,c)x)=[(a(b(cx)))p(b,c,x)p(a,bc,x)]\setminus
[(ab)(cx)p(ab,c,x)]$\par $p(u,x,(ux)\setminus (vx))[p(u,u\setminus
v,x)]^{-1}$, \\ because $u\setminus v=p(a,b,c)\in N(G)$ and
$p(u,u\setminus v,x)=e$.
\par Notice that $(2.1.1)$, $(2.1.2)$ and $(2.1.9)$ imply
$u\setminus (tu)=p$, where $t=t(a,b,c)$, $p=p(a,b,c)$, $u=a(bc)$ for
any $a$, $b$ and $c$ in $G$. Let $z\in G$, then there exists $x\in
G$ such that $z=ux$, that is $x=u\setminus z$. Therefore we deduce
that
\par $(2.8.6)$ $z\setminus (tz)=[x\setminus (px)]p(u,u\setminus
(tu),x)[p(u,x,(ux)\setminus (tux))]^{-1}$, \\ since $t\in N(G)$,
$p\in N(G)$, $(u\setminus (tu))x=(u\setminus (tux))[p(u,u\setminus
(tu),x)]^{-1}$ by $(2.2.2'')$; $~x\setminus (u\setminus
(tux))=[(ux)\setminus (tux))]p(u,x,(ux)\setminus (tux))$ by
$(2.2.2')$. Thus from Identities $(2.8.3)$, $(2.8.5)$ and $(2.8.6)$
it follows that a group $N_{0,0}=N_{0,0}(G)$ generated by $ \{
p(a,b,c), ~ t(a,b,c): a\in G, ~ b\in G, ~ c\in G \} $ satisfies
Condition $(2.7.1)$. From the joint continuity of the multiplication
and the mappings $Div_l$ and $Div_r$ it follows that the closure
$N_0$ of $N_{0,0}$ also satisfies $(2.7.1)$. Thus $N_0$ is a closed
normal subgroup in $G$. In view of Theorem 1.1 in Ch. IV, Section 1
in \cite{bruckb} a quotient loop $G/\cdot / N_0$ exists consisting
of all cosets $aN_0$, where $a\in G$.
\par  Then from
Conditions $(2.1.9)$, $(2.7.1)$ and $(2.7.2)$ it follows that for
each $a$, $b$, $c$ in $G$ the identities take place \par
$(aN_0)(bN_0)=(ab)N_0$ and
\par $((aN_0)(bN_0))(cN_0)=(aN_0)((bN_0)(cN_0))$ and
$eN_0=N_0$, \\ since $p(a,b,c)\in N_0$ and $t(a,b,c)\in N_0$ for all
$a$, $b$ and $c$ in $G$.
\par In view of Lemmas 2.2 and 2.3 $(aN_0)\setminus e= e/(aN_0)$,
consequently, for each $aN_0\in G/\cdot / N_0$ a unique inverse
$(aN_0)^{-1}$ exists. Thus the quotient $G/\cdot / N_0$ of $G$ by
$N_0$ is a group. Since the topology $\tau $ on $G$ is $T_1$ and
$N_0$ is closed in $G$, then the quotient topology $\tau _q$ on
$G/\cdot / N_0$ is also $T_1$. By virtue of Theorem 8.4 in
\cite{hew} this implies that $\tau _q$ is a $T_1\cap T_{3.5}$
topology on $G/\cdot / N_0$.
\par {\bf 2.9. Proposition.}  {\it Assume that $G$ is a $T_1$ topological
fan loop and functions $t$ and $p$ on $G$ are defined by Formulas
$(2.1.9)$. Then for each compact subset $S$ in $G$ and each open
neighborhood $V$ of $e$ there exists an open neighborhood $U$ of $e$
in $G$ such that
\par $(2.9.1)$ $t((u_1a)v_1,(u_2b)v_2,(u_3c)v_3)\in (Vt(a,b,c))\cap (t(a,b,c)V)$ and
\par $(2.9.2)$ $p((u_1a)v_1,(u_2b)v_2,(u_3c)v_3)\in (Vp(a,b,c))\cap (p(a,b,c)V)$
\\ for every $a$, $b$, $c$ in $S$ and $u_j$, $v_j$ in $\check{U}$ for each
$j\in \{ 1, 2, 3 \} $.}
\par {\bf Proof.} Take arbitrary fixed elements $f$, $g$ and $h$ in $S$.
From the joint continuity of the maps $t(a,b,c)$ and $p(a,b,c)$ in
the variables $a$, $b$ and $c$ in $G$ it follows that there exists
an open neighborhood $U_{f,g,h}$ of $e$ in $G$ and an open
neighborhood $W_{f,g,h}$ of $(f,g,h)\in S\times S\times S$ in
$G\times G\times G$ such that $(2.9.1)$ and $(2.9.2)$ are valid for
each $u_j$, $v_j$ in $\check{U}_{f,g,h}$, $j\in \{ 1, 2, 3 \} $, and
$(a,b,c)\in W_{f,g,h}$ (see Lemmas 2.4 and 2.6). Notice that
$S\times S\times S$ is compact in the Tychonoff product $G\times
G\times G$ of $G$ as the topological space (see Section 2.3 and
Theorem 3.2.4 in \cite{eng}). Hence an open covering $ \{ W_{f,g,h}:
~ f\in S, ~ g\in S, ~ h\in S \} $ of $S\times S\times S$ has a
finite subcovering $ \{ W_{f_i,g_i,h_i}: ~ i=1,...,n \} $, where $n$
is a natural number, $n\ge 1$. That is $S\times S\times S\subseteq
\bigcup_{i=1}^n W_{f_i,g_i,h_i}$. Then $\bigcap_{i=1}^n
U_{f_i,g_i,h_i}=:U$ is an open neighborhood of $e$ in $G$.
Therefore, Properties $(2.9.1)$ and $(2.9.2)$ are satisfied for
every $a$, $b$, $c$ in $S$ and $u_j$, $v_j$ in $\check{U}$ for each
$j\in \{ 1, 2, 3 \} $
\par We remind the following.

\par {\bf 2.10. Definition.} Let $G$ be a topological loop.
For a subset $U$ in $G$ it is put:
\par $(2.10.1)$ ${\cal L}_{U,G} := \{ (x,y)\in G\times G: ~ x\setminus y\in U \}
$ and
\par $(2.10.2)$ ${\cal R}_{U,G} := \{ (x,y)\in G\times G: ~ y/x \in U \}
$. \par A family of all subsets ${\cal L}_{U,G}$ (or ${\cal
R}_{U,G}$) with $U$ being an open neighborhood of $e$ will be
denoted by ${\cal L}_G$ (or ${\cal R}_G$ correspondingly).

\par {\bf 2.11. Proposition.} {\it Let $G$ be a $T_1$ topological locally compact fan
loop. Then a family ${\cal L}_G$ (or ${\cal R}_G$) induces a uniform
structure on $G$. A topology $\tau _1$ on $G$ provided by ${\cal
L}_G$ (or ${\cal R}_G$ respectively) is $T_1\cap T_{3.5}$ and
equivalent to the initial topology $\tau $ on $G$.}
\par {\bf Proof.} Let $(G, \tau )$ be a topological loop and let ${\cal
B}_e$ denote a base of its open neighborhoods at $e$. In view of
Lemma 2.5 ${\cal C}_l(U):= \{ xU: ~ x\in G \} $ is an open covering
of $G$ for each $U\in {\cal B}_e$. We put ${\cal C}_l^0 = \{ {\cal
C}_l(U): ~ U\in {\cal B}_e \} $ and ${\cal C}_l$ to be a family of
all coverings for each of which there exists a refinement of the
type ${\cal C}_l^0$. \par Below it is verified, that the family
${\cal C}_l$ satisfies Conditions $(UC1)$-$(UC4)$ of Section 8.1 in
\cite{eng}. If ${\cal A}\in {\cal C}_l$, ${\cal E}$ is a covering of
$G$ and ${\cal A}$ refines ${\cal E}$, then there exists $U\in {\cal
B}_e$ such that ${\cal C}_l(U)$ refines ${\cal A}$ and hence ${\cal
C}_l(U)$ refines ${\cal E}$. Thus $(UC1)$ is satisfied.
\par Let ${\cal A}_1$ and ${\cal A}_2$ belong to ${\cal C}_l$.
There are $U_1$ and $U_2$ in ${\cal B}_e$ such that ${\cal
C}_l(U_j)$ refines ${\cal A}_j$ for each $j\in \{ 1, 2 \} $. We put
$U=U_1\cap U_2$, consequently, $U\in {\cal B}_e$ and hence ${\cal
C}_l(U)$ refines both ${\cal C}_l(U_1)$ and ${\cal C}_l(U_2)$.
Therefore ${\cal C}_l(U)$ refines ${\cal A}_1$ and ${\cal A}_2$.
Thus $(UC2)$ also is satisfied.
\par Condition $(UC3)$ means that for each ${\cal A}\in {\cal C}_l$
there exists ${\cal E}\in {\cal C}_l$ such that ${\cal E}$ is a star
refinement of ${\cal A}$. In order to prove it, it evidently is
sufficient to prove that for each $U\in {\cal B}_e$ there exists
$U_1\in {\cal B}_e$ such that \par $(2.11.1)$ $St(xU_1,{\cal
C}_l(U_1))\subset xU$ for each $x\in G$, \\ where $St (M,{\cal A})$
denotes a star of a set $M$ with respect to ${\cal A}$ (see Section
5.1 in \cite{eng}).
\par Note that a map $f(x_1,x_2,x_3)=(x_1/x_2)x_3$ is the composition of
jointly continuous maps $G\times G\ni (x_1,x_2)\mapsto x_1/x_2\in G$
and $G\times G\ni (y,x_3)\mapsto yx_3\in G$, hence it is jointly
continuous from $G\times G\times G$ into $G$ and $f(e,e,e)=e$,
because $G$ is the topological loop (see Definition 2.1). The loop
$G$ is locally compact. Notice that for each open neighborhood $Q_1$
of $e$ in $G$ there exists an open neighborhood $Q_2$ of $e$ such
that its closure $cl_G(Q_2)$ is compact and $cl_G(Q_2)\subset Q_1$
by the corresponding Theorem 3.3.2 in \cite{eng} for topological
spaces. Hence for each open neighborhood $W$ of $e$ in $G$ there
exists an open neighborhood $U_0$ of $e$ in $G$ with the compact
closure $cl_G\check{U}_0$ such that $cl_G\check{U}_0$ is contained
in $W$ (see Lemma 2.6). \par Therefore for each $U\in {\cal B}_e$
there exists $V_1\in {\cal B}_e$ such that $f(V_1,V_1,V_1)\subset U$
and $cl_G(V_1)$ is compact. If for an arbitrary fixed element $x\in
G$ and some $x_1\in G$ the intersection $xV_1\cap x_1V_1\ne
\emptyset $ is non void, then there are $h_0$ and $h_1$ in $V_1$
such that $x_1=(xh_0)/h_1$. On the other hand, $x_1h\in x_1V_1$ for
each $h\in V_1$ and for each $y\in x_1V_1$ there exists $h\in V_1$
with $y=x_1h$, consequently, $x_1h=((xh_0)/h_1)h\in
((xV_1)/V_1)V_1$.
\par Using Identities $(2.2.3)$ and $(2.1.9)$ we get that
\par $(2.11.2)$ $x_1h=(x(h_0(e/h_1))p(x,h_0,e/h_1)$\par $p(e/h_1,h_1,h_1\setminus
e)[p((xh_0)/h_1,h_1,h_1\setminus e)]^{-1})h$. \par  We choose open
neighborhoods $V$ and $W$ of $e$ in $G$ such that
$\check{V}^2\subset W$ and $\check{W}^2\subset V_1$ by Lemma 2.6. In
view of the inclusion $(2.9.2)$ of Proposition 2.9 and Formula
$(2.11.2)$ there exists $U_1\in {\cal B}_e$ such that
$\check{U}_1\subset V$ and
\par $(2.11.3)$ $p((u_1a)v_1,(u_2b)v_2,(u_3c)v_3)\in (Vp(a,b,c))\cap (p(a,b,c)V)$
\\ for every $a$, $b$, $c$ in $cl_G(V_1)$ and $u_j$, $v_j$ in $\check{U}_1$ for each
$j\in \{ 1, 2, 3 \} $. This implies $(2.11.1)$ and hence $(UC3)$,
since $p(a,b,c)=e$ if either $a=e$ or $b=e$ or $c=e$.
\par It remains to prove that ${\cal C}_l$ also has the property
$(UC4)$. That is for each $x\ne y$ in $G$ there exists ${\cal A}\in
{\cal C}_l$ such that $\{ x, y \} \cap V\ne \{ x , y \} $ for each
$V\in {\cal A}$. It is sufficient to find an open neighborhood $U$
of $e$ in $G$ such that $x/U\cap y/U=\emptyset $, because this
implies $x_0U\cap \{ x, y \} \ne \{ x, y \} $ for each $x_0\in G$.
The loop $G$ is $T_1$. By virtue of Lemmas 2.5 and 2.6 and the joint
continuity of the multiplication and $Div_r$ in $G$ there is $U_1\in
{\cal B}_e$ such that $y\notin (xU_1)/U_1$, that is $xU_1\cap
yU_1=\emptyset $ by $(2.2.5)$. In view of Proposition 2.9 there
exists $U\in {\cal B}_e$ such that $(e/U)p(e/U,U,U\setminus
e)[p(a/U,U,U\setminus e)]^{-1}\subset U_1$ for each $a\in \{ x , y
\} $, since the two-point set $ \{ x, y \} $ is compact in $G$, for
each $W\in {\cal B}_e$ there exists $W_1\in {\cal B}_e$ such that
$e/W_1\subset W$. From $(2.2.3)$ it follows that $x/U\cap
y/U=\emptyset $. Therefore $\{ x, y \} \cap V\ne \{ x , y \} $  for
each $V\in {\cal C}_l(U)$.
\par  By virtue of Theorem 8.1.1 in \cite{eng} the uniformity ${\cal C}_l$ induces a $T_1$
topology $\tau _1$ on $G$. Note that the family ${\cal C}_l$
consists of open coverings of $G$ and that for each $x\in G$ and
each open neighborhood $V$ of $x$ in the initial topology $\tau $
there exists $U\in {\cal B}_e$ such that $xU\subset V$. Therefore
from the latter inclusion and $(2.11.1)$ it follows that the
topology $\tau _1$ induced by ${\cal C}_l$ coincides with the
initial topology $\tau $ on $G$. In view of Corollary 8.1.13 in
\cite{eng} $(G, \tau )$ is a Tychonoff space, that is $(G, \tau )$
is a completely regular space, $T_1\cap T_{3.5}$. Finally note that
${\cal C}_l^0={\cal L}_G$. Symmetrically the case ${\cal
C}_r^0={\cal R}_G$ is proved.

\par {\bf 2.12. Lemma.} {\it Suppose that $(G, \tau )$ is a $T_1$
topological loop, $S$ is a compact subset in $G$, $q$ is a fixed
element in $G$, $V$ is an open neighborhood of the unit element $e$.
Then there are elements $b_1,...,b_m$ in $G$ and an open
neighborhood $U$ of $e$ such that \par $(2.12.1)$ $\check{U}\subset
V$ and
\par $(2.12.2)$ $ \{ b_1\setminus (qU),...,b_m\setminus (qU) \} $ is
an open covering of $S$.}
\par {\bf Proof.} The multiplication is continuous on $G$, hence the left
shift mapping $L_b(x)=bx$ is continuous on $G$ in the variable $x$.
On the other hand, the mapping $Inv_l$ is continuous on $G$.
\par In view of $(2.1.1)$, $(2.1.2)$, Lemmas 2.5 and 2.6 and the
compactness of $S$ for each open neighborhood $U$ of $e$ in $G$ with
$\check{U}\subset V$ there are $b_1,...,b_m$ in $G$ such that $\{
b_1\setminus (qU),...,b_m\setminus (qU) \} $ is an open covering of
$S$.

\par {\bf 2.13. Corollary.} {\it Let $G$ be a $T_1$ topological loop. Then for each open neighborhood
$W$ of $e$ in $G$ there exists an open neighborhood $U$ of $e$ such
that $\check{U}\subset W$ and \par $(2.13.1)$ $(\forall x$ $\forall
y$ $((x\in G) \& (y\in G)\& (x\setminus y \in U))) \Rightarrow (y\in
xW)$ and
\par $(2.13.2)$ $(\forall x$ $\forall y$ $((x\in G) \&
(y\in G)\& (y/x \in U))) \Rightarrow (y\in Wx)$.}
\par {\bf Proof.} This follows from Lemmas 2.6 and 2.12, $(2.1.1)$ and $(2.1.2)$.

\par {\bf 2.14. Theorem.} {\it Let $G$ and $H$ be $T_1$
topological fan  loops (see Definition 2.1) and let $f: G\to H$ be a
continuous map so that for each open neighborhood $V$ of a unit
element $e_H$ in $H$ a compact subset $K_V$ in $G$ exists such that
$f(G- K_V)\subset V$. Then $f$ is uniformly $({\cal L}_G, {\cal
L}_H)$ continuous and uniformly $({\cal R}_G, {\cal R}_H)$
continuous (see also Definition 2.10).}
\par {\bf Proof.} Since the multiplication in $H$ is continuous, then for each open
neighborhood $Y$ of $e_H$ there exists an open neighborhood $X$ of
$e_H$ such that $X^2\subset Y$. In view of Lemma 2.6 there exists an
open neighborhood $V_1$ of $e_H$ in $H$ such that
$\check{V}_1^2\subset V$, where $A^2=AA$ for a subset $A$ in $H$. By
the conditions of this theorem a compact subset $K_{V_1}$ in $G$
exists such that $f(G- K_{V_1})\subset V_1$.
\par For a subset $A$ of the loop $G$ let \par $(2.14.1)$
$P(A)=(P_0(A)\cup \{ e \} )(P_0(A)\cup \{ e \} )$, \par where
$P_0(A)=A\cup Inv_l(A)\cup Inv_r(A)$, \\ hence $A\subset P_0(A)$ and
$P_0(A)\cup \{ e \} \subset P(A)$. Then $S_1=P(K_{V_1})$ is a
compact subset in $G$, since the mappings $Inv_l$ and $Inv_r$ are
continuous on $G$ and the multiplication is jointly continuous on
$G\times G$ (see Theorems 3.1.10, 8.3.13-8.3.15 in \cite{eng}),
hence $R_1=P(f(S_1))$ is compact in $H$.
\par  By virtue of Proposition 2.9 there
exists an open neighborhood ${V'}_2$ of $e_H$ in $H$ such that
\par $(2.14.2)$ $[t_H((V_2a)V_2,(V_2b)V_2,(V_2c)V_2)V_2]\cup
[V_2t_H((V_2a)V_2,(V_2b)V_2,(V_2c)V_2)]$\par $\subset
(V_3t_H(a,b,c))\cap (t_H(a,b,c)V_3)$ and
\par  $[p_H((V_2a)V_2,(V_2b)V_2,(V_2c)V_2)V_2]\cup
[V_2p_H((V_2a)V_2,(V_2b)V_2,(V_2c)V_2)]$\par $\subset
(V_3p_H(a,b,c))\cap (p_H(a,b,c)V_3)$
\\ for every $a$, $b$, $c$ in $R_1$, where $\check{V}_3^2\subset V_1$, $ ~ V_2=\check{V'}_2$,
$V_3$ is an open neighborhood of $e$ in $H$. For $V_2$ there exists
a compact subset $K_{V_2}$ in $G$ such that $f(G- K_{V_2})\subset
V_2$ by the conditions of this theorem. If $A$ and $B$ are compact
subsets in $G$, then their union $A\cup B$ is also compact.
Therefore it is possible to choose $K_{V_2}$ such that
$K_{V_1}\subset K_{V_2}$, since $V_2\subset V_1$ and $(G- A)- B=G-
(A\cup B)\subset G- A$. We take $S_2=P(K_{V_2})$ by Formula
$(2.14.1)$, consequently, $S_1\subset S_2$, since $K_{V_1}\subset
K_{V_2}$.
\par From the continuity of the map $f$ and Lemmas 2.5, 2.6
it follows that for each $x\in G$ open neighborhoods $W_{x,l}$ and
$W_{x,r}$ of $e$ in $G$ exist such that
\par $f(x\check{W}_{x,l}^2)\subset (f(x)V_2)$ and $f(\check{W}_{x,r}^2x)\subset
(V_2f(x))$, consequently,
\par $(2.14.3)$ $f(x\check{W}_x^2)\subset (f(x)V_2)$ and
$f(\check{W}_x^2x)\subset (V_2f(x))$ \\ for an open neighborhood
$W_x=W_{x,l}\cap W_{x,r}$ of $e$ in $G$. The compactness of $S_2$
imply that coverings $ \{ xW_x: ~ x\in S_2 \} $ and $ \{ W_yy: ~
y\in S_2 \} $ of $S_2$ have finite subcoverings $ \{ x_jW_{x_j}: ~
x_j\in S_2, ~ j=1,...,n \} $ and $ \{ W_{y_i}y_i: ~ y_i\in S_2, ~
i=1,...,m \} $. Hence \par $(2.14.4)$ $W=\bigcap_{j=1}^nW_{x_j}\cap
\bigcap_{i=1}^m W_{y_i}$ \\ is an open neighborhood of $e$ in $G$.
Therefore according to Proposition 2.9 there exists an open
neighborhood $U'$ of the unit element $e$ in $G$ such that
\par $(2.14.5)$ $[t_G((Ua)U,(Ub)U,(Uc)U)U]\cup [Ut_G((Ua)U,(Ub)U,(Uc)U)]
$\par $\subset [W_3t_G(a,b,c)]\cap [t_G(a,b,c)W_3]$ and
\par $[p_G((Ua)U,(Ub)U,(Uc)U)U]\cup [Up_G((Ua)U,(Ub)U,(Uc)U)]
$\par $\subset [W_3p_G(a,b,c)]\cap [p_G(a,b,c)W_3]$ \\ for every
$a$, $b$, $c$ in $S_2$, where $U=\check{U'}$, $~ W_0$ and $W_3$ are
open neighborhoods of $e$ in $G$ such that $\check{W}_3^2\subset
W_0$ and $\check{W}_0^2\subset W$.
\par Let now $x$ and $y$ in $G$ be such that $x\setminus y\in U$.
Then Formula $(2.2.4)$ imply that
\par $(2.14.6)$ $y\in xU$.
\par There are several options. Consider at first the
case $x\in K_{V_2}$. From Formulas $(2.14.4)$-$(2.14.6)$ and
Corollary $(2.13)$ it follows that there exists $j\in \{ 1,...,n \}$
such that $x\in x_jW_{x_j}$ and $y\in x_jW_{x_j}^2$. Therefore,
Formulas $(2.14.2)$ and $(2.14.3)$ imply that $f(x)\setminus f(y)\in
V$.
\par From $x\setminus y\in U$ and Identities $(2.2.4)$ it follows that
$y=xu$ for a unique $u\in U$. Hence \par $(2.14.7)$ $x=[t(y,e/u,u)]^{-1} y(e/u)$ \\
according to Identities $(2.2.3)$, $(2.2.5)$.
\par If $y\in K_{V_2}$, then similarly from Formulas $(2.14.4)$-
$(2.14.7)$ and Corollary $(2.13)$ it follows that there exists $k\in
\{ 1,...,n \} $ such that $y\in x_kW_{x_k}$ and $x \in
x_kW_{x_k}^2$, since $t(a,b,e)=t(a,e,b)=t(e,a,b)=e$ for each $a$ and
$b$ in $G$. Therefore, $f(x)\setminus f(y)\in V$ by Formulas
$(2.14.2)$ and $(2.14.3)$, since $S_2=P(K_{V_2})$ (see Formula
$(2.14.1)$).
\par It remains the case $x\in G- K_{V_2}$ and $y\in
G- K_{V_2}$. Therefore $f(x)\in V_2$ and $f(y)\in V_2$. According to
the choice of $R_1$ we have $e_H\in R_1$. From Condition $(2.14.2)$,
Identity $(2.2.4)$ and the inclusion $\check{V}_1^2\subset V$, it
follows that $f(x)\setminus f(y)\in V$. Taking into account the
inclusion $K_{V_1}\subset K_{V_2}$ we get that $f$ is uniformly
$({\cal L}_G, {\cal L}_H)$ continuous.
\par The uniform $({\cal R}_G, {\cal R}_H)$ continuity is proved
analogously using the finite subcovering $ \{ W_{y_i}y_i: ~ y_i\in
S_2, ~ i=1,...,m \} $ and Corollary 2.13.

\par {\bf 2.15. Corollary.} {\it Let $G$ be a $T_1$ topological locally compact fan
loop and let $f\in C_0(G)$ and let $H=({\bf C}, +)$ be the complex
field $\bf C$ considered as an additive group. Then $f$ is uniformly
$({\cal L}_G, {\cal L}_H)$ continuous and uniformly $({\cal R}_G,
{\cal R}_H)$ continuous.}

\section{Left invariant measures.}
\par {\bf 3.1. Notation.} For a completely regular topological space $X$ by $C_b(X)$
is denoted the Banach space of all continuous bounded functions $f$
from $X$ into the complex field $\bf C$ supplied with the norm
\par $(3.1.1)$ $\| f \|_X = \sup_{x\in X} |f(x)|<\infty $.
\par We put
\par $(3.1.2)$ $C_0(X) := \{ f \in C_b(X): ~ \forall ~ \epsilon
>0, ~ \exists ~ S\subset X, ~ S \mbox{ is compact, }$\par $ ~ \forall ~ x\in
X- S, ~ |f(x)|<\epsilon \} $ and
\par $(3.1.3)$ $C_{0,0}(X) := \{ f \in C_b(X):
~ \exists S\subset X, ~ S \mbox{ is compact, }$\par $ ~ \forall x\in
X- S, ~ f(x)=0 \} $ and
\par $(3.1.4)$ $C_{0,0}^+(X) = \{ f\in C_{0,0}(X): ~ \forall x\in X, ~ f(x)\ge 0 \}
$.
\par Let $G$ be a loop. For a function $f: G\to \bf C$ and an element $b\in G$
let $L_bf(x)=\mbox{}_bf(x)=f(bx)$ and $R_bf(x)=f_b(x)=f(xb)$ for
each $x\in G$. Consider a support $S_f := cl_G \{ x\in G: ~ f(x)\ne
0 \} $ of $f\in C_b(G)$, where $cl_G(A)$ denotes the closure of a
subset $A$ in $G$.

\par {\bf 3.2. Lemma.} {\it Let $(G, \tau )$ be a
$T_1$ topological locally compact fan  loop, let also $f$ and $\phi
$ belong to $C_{0,0}^+(G)$ and $\phi $ be not identically zero (see
Notation 3.1). Then there exist a natural number $m>0$, elements
$b_1,...,b_m$ in $G$ and positive constants $c_1,...,c_m$ such that
\par $(3.2.1)$ $\forall x\in G$, $ ~f(x)\le \sum_{j=1}^m c_j L_{b_j}\phi
(x).$}
\par {\bf Proof.} Since $f\in C_{0,0}^+(G)$, then the support $S_f$ is compact.
The function $\phi $ is not null, hence there exists $q\in G$ such
that $\phi (q)>0$. From Lemma 2.5 and from the continuity of the
function $\phi $ it follows that there exists an open neighborhood
$qV$ of $q$ such that $\phi (x)> \phi (q)/2$ for reach $x\in qV$,
where $V$ is an open neighborhood of the unit element $e$. By virtue
of Lemma 2.12 there exists an open neighborhood $U$ of $e$ and
elements $b_1,...,b_m$ in $G$ such that $\check{U}\subset V$ and for
each $x\in S_f$ there exists $j\in \{ 1,...,m \} $ such that $x\in
b_j\setminus (qU)$.
\par Therefore, $$f(x)\le \| f \| _G (2/\phi (q))\sum_{j=1}^m \phi
(b_jx)$$ for each $x\in G$ according to $(2.2.4)$, so it is
sufficient to take $c_j\ge \| f \| _G (2/\phi (q))$ for each
$j=1,...,m$.

\par {\bf 3.3. Corollary.} {\it Let the conditions of Lemma 3.2 be
satisfied and let
$$(3.3.1) \quad (f:\phi ) := \inf \{ \sum_{j=1}^m c_j: ~ \exists
~ \{ b_1,...,b_m \} \subset G, ~ \exists ~ \{ c_1,...,c_m \} \subset
(0, \infty ),$$ $$ ~ \forall ~ x\in G, ~ f(x)\le \sum_{j=1}^m c_j
L_{b_j}\phi (x) \} .$$ Then $(f: \phi )\le 2m\| f \| _G /\phi (q)$
in the notation of Lemma 3.2.}

\par \par {\bf 3.4. Lemma.} {\it Assume that the conditions of Lemma 3.2 are
fulfilled, then for each $b\in G$ \par $(3.4.1)$ $(\mbox{}_b f: \phi
)=(f: \phi ^b)$;
\par $(3.4.2)$ $(f: \mbox{}_b\phi
)=(f^b: \phi )$,
\\ where $f^b(x)=f(b\setminus
x)$ for each $x\in G$; particularly, \par $(3.4.1')$
$(\mbox{}_{\gamma } f: \phi )=(f: \phi )$ and
\par $(3.4.2')$
$(f: \mbox{}_{\gamma } \phi )=(f: \phi )$ for each $\gamma \in
N(G)$;
\par $(3.4.3)$ $(\alpha f: \phi )=\alpha (f: \phi )$ for each
$\alpha \ge 0$;
\par $(3.4.4)$
$((f_1+f_2): \phi )\le (f_1: \phi ) +(f_2: \phi )$ for every $f_1$
and $f_2$ in $C_{0,0}^+(G)$.
\par $(3.4.5)$ If $f(x)\le f_1(x)$ for each $x\in G$, then $(f: \phi
)\le (f_1: \phi )$.}
\par {\bf Proof.} Let $c_1,...,c_m$ in $(0, \infty )$ and
$b_1,...,b_ m$ in $G$ be such that $$(3.4.6)\quad \mbox{}_bf(x)\le
\sum_{j=1}^mc_jL_{b_j}\phi (x)$$ for each $x\in G$. From Formulas
$(2.2.4)$ and $(3.4.6)$ by changing of a variable $y=bx$ it follows
that
$$(3.4.7)\quad f(y)\le
\sum_{j=1}^mc_jL_{b_j}\phi (b\setminus y)$$ for each $y\in G$. From
$(3.4.7)$ it follows $(3.4.1)$. Similarly from the inequality
$$(3.4.8)\quad f(x)\le
\sum_{j=1}^mc_jL_{b_j}(L_b\phi (x))$$ for each $x\in G$ we infer
that
$$(3.4.9)\quad f(b\setminus y)\le \sum_{j=1}^mc_jL_{b_j}\phi (y)$$ for each $y\in
G$. Thus $(3.4.9)$ implies Equality $(3.4.2)$.
\par In particular, if $\gamma \in N(G)$, then
$b_j(\gamma \setminus y)=(b_j\gamma ^{-1})y$ and $b_j(\gamma
y)=(b_j\gamma )y$ for each $y$ and $b_j$ in $G$ by Condition
$(2.1.8)$ and Formulas $(2.2.2)$ and $(2.3.1)$. Hence $(3.4.7)$
transforms into to
$$f(y)\le \sum_{j=1}^mc_jL_{b_j\gamma ^{-1}}\phi (y)$$ and $(3.4.8)$
into
$$f(x)\le \sum_{j=1}^mc_jL_{b_j\gamma }\phi (x)$$ with $\gamma \in N(G)$
instead of $b$. This implies Equalities $(3.4.1')$, $(3.4.2')$.
\par Properties $(3.4.3)$ and $(3.4.4)$ evidently follow from
Formula $(3.3.1)$.
\par For proving Property $(3.4.5)$ note that if $f(x)\le f_1(x)$ for each $x\in
G$, then from $f_1(x)\le \sum_{j=1}^mc_jL_{b_j}\phi (x)$ for each
$x\in G$ it follows that $f(x)\le \sum_{j=1}^mc_jL_{b_j}\phi (x)$
for each $x\in G$, consequently, $(f: \phi )\le (f_1: \phi )$.

\par {\bf 3.5. Notation.} Let $\phi $, $f_0$ and $f$ belong to
$C_{0,0}^+(G)$ and $\phi $ and $f_0$ be not null, where $G$ is a
$T_1$ topological locally compact fan loop. We consider a functional
$$(3.5.1)\quad J_{\phi ,f_0}(f):= \frac{(f: \phi )}{(f_0: \phi )}.$$
\par Assume that \par $(3.5.2)$ there exists a compact subgroup $N_0=N_0(G)$ in $N(G)$
such that \par $t(a,b,c)\in N_0$ and $p(a,b,c)\in N_0$ for every
$a$, $b$ and $c$ in $G$.
\par Then we denote by $\Upsilon (G,N_0)$ a family of all non null functions $h$ in $C_{0,0}^+(G)$ such that
\par $(3.5.3)$  $h(\gamma a)=h(a)$
for each $a\in G$ and $\gamma \in N_0$.
\par Evidently Condition $(3.5.3)$ for $h\in C_{0,0}^+(G)$ is equivalent to
\par $(3.5.4)$ $h(a\gamma )=h(a)$ for each $a\in G$ and $\gamma \in N_0$, \\ since
$aN_0=N_0a$ for each $a\in G$ according to Theorem 2.8.

\par {\bf 3.6. Lemma.} {\it Let $G$ be a $T_1$ topological locally compact fan  loop
satisfying Condition $(3.5.2)$, $f$ and $\phi $ be in $C^+_{0,0}(G)$
and $\omega \in \Upsilon (G,N_0)$ (see Condition $(3.5.3)$), $\phi $
be non null. Then
\par $(3.6.1)$ $(f: \phi ) \le (f: \omega )
(\omega : \phi )$.}
\par {\bf Proof.} If $b$ is a fixed element in $G$ and there are elements
$b_1,...,b_m$ in $G$ and positive constants $c_1,...,c_m$ such that
$$(3.6.2)\quad \mbox{}_b\omega(x)\le \sum_{j=1}^mc_j\phi (b_jx)$$ for
each $x\in G$, then
$$(3.6.3)\quad \mbox{}_b\omega(x)\le
\sum_{j=1}^mc_j\phi (b_jx\gamma )$$ for each $x\in G$ and $\gamma
\in N_0$, since $N_0\subset N(G)$ and $\mbox{}_b\omega(x\gamma
)=\mbox{}_b\omega(x)$ for each $x\in G$ and $\gamma \in N_0$ by
$(3.5.4)$ equivalent to $(3.5.3)$.
\par By the conditions of this lemma $N_0$ is a compact group.
Therefore there exists a Haar measure $\lambda $ on the Borel
$\sigma $-algebra ${\cal B}(N_0)$ of $N_0$ and with values in the
unit segment $[0,1]$ such that $\lambda (N_0)=1$, $\lambda
(sA)=\lambda (A)$ and $\lambda (As)=\lambda (A)$ for each $s\in N_0$
and $A\in {\cal B}(N_0)$ (see Theorems 15.5, 15.9 and 15.13 and
Subsection 15.8 in \cite{hew}). In view of this, Conditions
$(3.1.3)$ and $(3.1.4)$ and Corollary 2.15 a function
$$(3.6.4)\quad \phi ^{[\lambda ]}(x):=\int_{N_0} \phi (\gamma x)\lambda (d\gamma )$$
on $G$ is nonzero and belongs to $C_{0,0}^+(G)$, since $N_0S_{\phi
}$ is a compact subset in $G$ by Lemma 2.6, where $S_{\phi }$ is a
compact support of $\phi $. From Formula $(3.6.4)$ it follows that
\par $(3.6.5)$ $\phi ^{[\lambda ]}(\beta
x)=\phi ^{[\lambda ]}(x)$ for each $\beta \in N_0$ and $x\in G$, \\
since the measure $\lambda $ is left and right invariant $\lambda
(\beta A)=\lambda (A)= \lambda (A\beta )$ for each $\beta \in N_0$
and each Borel subset $A$ in $N_0$. Hence $\phi ^{[\lambda ]}\in
\Upsilon (G,N_0)$, since $S_{\phi }N_0$ is compact, and since
Conditions $(3.5.3)$ and $(3.5.4)$ are equivalent, where $S_{\phi }$
is the support of $\phi $ (see Subsection 3.2). From $(3.6.4)$,
$(3.6.5)$, $(3.5.3)$, $(3.5.4)$ and the Fubini theorem it follows
that
$$(3.6.4')\quad \phi ^{[\lambda ]}(x)=\int_{N_0} \phi (x\beta )\lambda (d\beta )\mbox{,
since}$$
$$\phi ^{[\lambda ]}(x)=\int_{N_0}(\int_{N_0} \phi (\gamma x\beta )\lambda (d\gamma
))\lambda (d\beta )$$
$$=\int_{N_0}(\int_{N_0} \phi (\gamma x\beta )\lambda (d\beta ))\lambda (d\gamma
)=\int_{N_0} \phi (x\beta )\lambda (d\beta )),$$ because $\int_{N_0}
\phi (x\gamma \beta )\lambda (d\beta ) =\int_{N_0} \phi (x\beta
)\lambda (d\beta )$ for each $\gamma \in N_0(G)$.
\par Integrating both sides of Inequality $(3.6.3)$ and utilizing
Formulas $(3.6.4)$, $(3.6.4')$ we infer that
$$(3.6.6)\quad \mbox{}_b\omega(x)\le
\sum_{j=1}^mc_j\phi ^{[\lambda ]}(b_jx)$$ for each $x\in G$. On the
other hand,
$$\int_{N_0} (\sum_{j=1}^mc_j\mbox{ }_{b_j}\phi )(x\gamma )\lambda
(d\gamma )=(\sum_{j=1}^mc_j\mbox{ }_{b_j}\phi )^{[\lambda ]
}(x)=\sum_{j=1}^mc_j\mbox{ }_{b_j}\phi ^{[\lambda ]}(x),$$ hence for
each $x\in G$ there exists $\gamma \in N_0$ such that
$$(\sum_{j=1}^mc_j\mbox{ }_{b_j}\phi )(x\gamma )\ge
\sum_{j=1}^mc_j\mbox{ }_{b_j}\phi ^{[\lambda ]}(x).$$ Thus vice
versa from $\omega \in \Upsilon (G,N_0)$ and $(3.6.6)$ it follows
$(3.6.3)$ and hence $(3.6.2)$, consequently,
\par $(3.6.7)$ $(\mbox{}_b\omega: \phi ^{[\lambda ]})= (\mbox{}_b\omega: \phi )$.
\par Let $a_1,...,a_n$ in $G$ and positive constants
$q_1,...,q_n$ be such that $$(3.6.8)\quad \mbox{}_b\omega(x)\le
\sum_{j=1}^nq_j\phi ^{[\lambda ]}(a_jx)$$ for each $x\in G$ (see
Lemma 3.2). From Formulas $(2.2.2)$, $(3.6.5)$, $(3.6.8)$ and
Conditions $(3.5.2)$, $(3.5.3)$, $(3.5.4)$ we deduce that
$$(3.6.9)\quad \omega(y)\le
\sum_{j=1}^nq_j\phi ^{[\lambda ]}((a_j(b\setminus e))y
[p(a_j,b\setminus e,y)]^{-1} p(b,b\setminus e, y) )$$
$$=\sum_{j=1}^nq_j\phi ^{[\lambda ]} (d_jy)$$ for each $y\in G$,
where $d_j=a_j(b\setminus e)$ for each $j$. Therefore
\par $(\mbox{}_b\omega: \phi ^{[\lambda ]})\le (\omega: \phi ^{[\lambda ]})$
for each $b\in G$. Notice that \par $(3.6.10)$ $L_cL_{c\setminus
e}\omega(x)=\omega(x)$ for each $c$ and $x$ in $G$ by Lemmas 2.2,
2.3 and Condition $(3.5.3)$. Therefore we analogously get \par
$(\omega: \phi ^{[\lambda ]}) \le (\mbox{}_c\omega : \phi ^{[\lambda
]})$ for each $c\in G$. Thus
\par $(3.6.11)$ $(\mbox{}_b\omega: \phi ^{[\lambda ]})= (\omega: \phi
^{[\lambda ]})$ for each $b\in G$. \par From $(3.6.7)$ and
$(3.6.11)$ it follows that
\par $(3.6.12)$ $(\mbox{}_b\omega: \phi )= (\omega: \phi )$ for each $b\in
G$. \par If $c_1,...,c_n$, $h_1,...,h_k$ in $(0,\infty )$ and
$a_1,..,a_k$, $g_1,...,g_n$ in $G$ are such that
$$(3.6.13)\quad f(x)\le \sum_{j=1}^kh_jL_{a_j}\omega (x) \mbox{  and}$$
$$(3.6.14)\quad \omega (x)\le \sum_{i=1}^n c_i L_{g_i}\phi (x)$$
for each $x\in G$ (see Lemma 3.2). Then from $(3.5.3)$, $(3.6.7)$,
$(3.6.12)$-$(3.6.14)$ and Lemma 2.2 we infer that
$$(3.6.15)\quad f(x)\le \sum_{j=1}^k h_j \sum_{i=1}^n c_i
L_{g_i}L_{a_j}\phi (x)= \sum_{j=1}^k h_j \sum_{i=1}^n c_i \phi
((g_ia_j)x).$$  Apparently $(3.6.15)$ implies $(3.6.1)$.

\par {\bf 3.7. Lemma.} {\it Let $G$ be a
$T_1$ topological locally compact fan loop, $\phi $, $f_0$ be
nonzero functions belonging to $C_{0,0}^+(G)$. Then for each
functions $f$, $f_1$ in $C_{0,0}^+(G)$ and $\alpha \ge 0$
\par $(3.7.1)$ $J_{\phi ,f_0}(\alpha f)=\alpha J_{\phi ,f_0}(f)$ and
\par $(3.7.2)$ $J_{\phi ,f_0}(f+f_1)\le J_{\phi ,f_0}(f) + J_{\phi
,f_0}(f_1)$ and
\par $(3.7.3)$ if $f(x)\le f_1(x)$ for each $x\in G$, then
$J_{\phi ,f_0}(f)\le J_{\phi ,f_0}(f_1)$.
\par Moreover, if $G$ satisfies Condition $(3.5.2)$ and $f_0\in \Upsilon (G,N_0)$
(see Condition $(3.5.3)$), then \par $(3.7.4)\quad  J_{\phi ,f_0}(f)
\le (f: f_0).$ }
\par {\bf Proof.} Properties $(3.7.1)$ and $(3.7.2)$ follow
immediately from $(3.4.3)$ and $(3.4.4)$. Property $(3.7.3)$ follows
from Property $(3.4.5)$.
\par  Applying Inequality $(3.6.1)$ and Formula
$(3.5.1)$ we infer Inequality $(3.7.4)$, since $J_{\phi ,
f_0}(f_0)=1$.

\par {\bf 3.8. Lemma.} {\it Assume that $G$ is a
$T_1$ topological locally compact fan loop, functions $\phi $, $f_0$
and $f$ belong to $C_{0,0}^+(G)$ and $\phi $ and $f_0$ are not null.
Then mappings $J_{\phi ,f_0}(\mbox{}_bf)$ and $J_{\phi ,f_0}(f_b)$
are continuous in the variable $b$ in $G$.}
\par {\bf Proof.} For each $x$, $b_1$ and $b_2$ in $G$ we have
$\mbox{}_{b_1}f(x)-\mbox{}_{b_2}f(x)=f(b_1x)-f(b_2x)$. In view of
Corollary 2.15 for each $\epsilon >0$ there exists an open of the
form $(2.6.1)$ neighborhood $U$ of $e$ in $G$ with a compact closure
$cl_G(U)$ for which \par $(3.8.1)$ $|f(b_1x)-f(b_2x)|<\epsilon $ for
each $x$, $b_1$ and $b_2$ in $G$ such that $(b_2x)\setminus
(b_1x)\in U$.
\par On the other hand, a support $S_f$ of $f$ is compact, consequently,
$bS_f=L_bS_f$ is compact for each $b\in G$. Let $b_1$ be fixed. For
each $x\in G$ with $b_1x\in S_f$ there exists an open of the form
$(2.6.1)$ neighborhood $W_x$ of $e$ in $G$ such that
$(b_2x)\setminus (b_1x)\in U$ for each $b_2x\in (b_1W_x)x\cap
b_1(xW_x)$ according to Lemmas 2.2, 2.4, 2.5, Proposition 2.9 and
Formula $(2.14.3)$. For an open covering $\{ (b_1W_x)x\cap
b_1(xW_x): ~ b_1x\in S_f, ~ x\in G \} $ of $S_f$ there exists a
finite subcovering $\{ (b_1W_{x_j})x_j\cap b_1(x_jW_{x_j}): ~
b_1x_j\in S_f, ~ x_j\in G, ~ j=1,...,m \} $ (see also Lemma 2.5),
since the subset $S_f$ is compact.
\par We take $W_0=U\cap \bigcap_{j=1}^m W_{x_j}$ and choose an open
of the form $(2.6.1)$ neighborhood $W$ of $e$ in $G$ with a compact
closure $cl_G(W)$ contained in $W_0$ (see Theorem 3.3.2 in
\cite{eng} and Formula $(2.14.3)$), because $G$ is locally compact.
\par In view of Proposition 2.9 and Lemma 2.6 there exists an open
neighborhood $V'$ of $e$ in $G$ with $V=\check{V'}$ and a compact
closure $cl_G(V)$ such that
\par $(3.8.2)$ $[t((Va)V,(Vb)V,(Vc)V)V]\cup [Vt((Va)V,(Vb)V,(Vc)V)]
$\par $\subset [t(a,b,c)W_1]\cap [W_1t(a,b,c)]$ and
\par $[p((Va)V,(Vb)V,(Vc)V)V]\cup [Vp((Va)V,(Vb)V,(Vc)V)]
$\par $\subset [p(a,b,c)W_1]\cap [W_1p(a,b,c)]$ \\
for each $a$, $b$ and $c$ in $S$, where $\check{W}_1^2\subset W$, $
~ W_1$ is an open neighborhood of $e$ in $G$,
\par $S=P(S_1)$, $S_1=S_2 \cup cl_G(U)$, \par $S_2=\{ y\in G: ~
y=(b_1u)x, ~ u\in cl_G(U), b_1x\in S_f \} $ \\ (see Formula
$(2.14.1)$), since $S$ is compact, $t(a,b,c)=e$ and $p(a,b,c)=e$ if
$e\in \{ a, b, c \} $. For $b_1x\notin S_f$ and $b_2x\notin S_f$
certainly $f(b_1x)-f(b_2x)=0$. So remain two cases either $b_1x\in
S_f$ or $b_2x\in S_f$ which are similar to each other up to a
notation. From Formulas $(2.2.5)$ it follows that $b_2x\in (b_1V)x$
is equivalent to $b_2\in b_1V$. Hence Lemma 2.2 and Inclusion
$(3.8.2)$ provide that $(b_2x)\setminus (b_1x)\in U$ for each
$b_2\in b_1V$ and $b_1x\in S_f$. \par Let $w\in C_{0,0}^+(G)$ be a
function such that $w(y)=1$ for each $y\in (cl_G(U)S_f)cl_G(U)$.
Then we deduce that $|f(b_1x)-f(b_2x)|<\epsilon w(x)$ for each $x$,
$b_1$ and $b_2$ in $G$ such that $b_2\in b_1V$ and with $b_1x\in
S_f$. \par Therefore for each $\epsilon
>0$ there exists and open neighborhood $V$ of $e$ in $G$
such that $|(\mbox{}_{b_1}f: \phi )- (\mbox{}_{b_2}f: \phi
)|<\epsilon (w:\phi )$ for each $b_2\in b_1V$, \\ consequently,
\par $(3.8.3)$ $|J_{\phi ,f_0}(\mbox{}_{b_1}f)- J_{\phi
,f_0}(\mbox{}_{b_2}f)|<\epsilon J_{\phi ,f_0}(w)$ \\ according to
Formula $(3.5.1)$, since $(f_0: \phi )>0$. Thus the mapping $J_{\phi
,f_0}(\mbox{}_{b}f)$ is continuous in the parameter $b\in G$, since
$0<J_{\phi ,f_0}(w)<\infty $ (see Lemmas 3.2, 3.7 and Corollary
3.3).
\par The case $J_{\phi ,f_0}(f_b)$ is proved symmetrically.

\par {\bf 3.9. Theorem.} {\it Assume that $G$ is a
$T_1$ topological locally compact fan loop satisfying Condition
$(3.5.2)$, $\phi $, $f$ and $f_1$ are nonzero functions belonging to
$C^+_{0,0}(G)$ and $f_0\in \Upsilon (G,N_0)$ (see $(3.5.3)$), then
\par $(3.9.1)\quad  (f_0:f ) ^{-1}\le J_{\phi ,f_0}(f) \le (f:f_0)$ and
\par $(3.9.2)\quad  (f_1: f_0)^{-1}(f_0: f ) ^{-1}\le J_{\phi ,f_1}(f)
\le (f: f_0) (f_0: f_1)$.}
\par {\bf Proof.} The right
inequality in $(3.9.1)$ follows from the inequality $(3.7.4)$.
\par Formulas $(3.6.7)$ and $(3.6.12)$ imply that
\par $(3.9.3)$ $(\mbox{}_bf_0: f)=
(f_0: f^{[\lambda ]})$ and $(\mbox{}_bf^{[\lambda ]}: \phi )=
(f^{[\lambda ]}: \phi ^{[\lambda ]})$ for each $b\in G$.
\par Let $c_1,...,c_k$, $h_1,...,h_n$ in
$(0,\infty )$ and $a_1,..,a_k$, $g_1,...,g_n$ in $G$ be such that
$$(3.9.4)\quad f_0(x)\le \sum_{j=1}^kc_jf^{[\lambda ]}(a_jx) \mbox{  and}$$
$$(3.9.5)\quad f^{[\lambda ]}(x)\le \sum_{i=1}^n h_i \phi ^{[\lambda ]}(g_ix)$$
for each $x\in G$ (see Lemma 3.2). Then from Identity $(2.1.9)$,
Inequalities $(3.9.4)$, $(3.9.5)$ and Conditions $(3.5.3)$,
$(3.5.4)$ we deduce that
$$(3.9.6)\quad f_0(x)\le \sum_{j=1}^k c_j \sum_{i=1}^n h_i
\phi ^{[\lambda ]}((g_ia_j)x[p(g_i,a_j,x)]^{-1})=\sum_{j=1}^k c_j
\sum_{i=1}^n h_i \phi ^{[\lambda ]}((g_ia_j)x).$$ Suppose that there
are $y_1,...,y_k\in G$ and $q_1,...,q_k\in ( 0, \infty )$ such that
$$(3.9.7)\quad f(x)\le \sum_{i=1}^k q_i \phi (y_ix)$$ for each
$x\in G$. Taking the integral $\int_{N_0}f(x\gamma )\lambda (d\gamma
)$ and similarly for the right side (see Formulas $(3.6.4)$ and
$(3.6.4')$), we get from Inequality $(3.9.7)$ that
$$f^{[\lambda ]}(x)\le \sum_{i=1}^k q_i
\phi ^{[\lambda ]}(y_ix)$$ for each $x\in G$ (see Lemma 3.2). Hence
\par $(3.9.8)$ $(f^{[\lambda ]}: \phi ^{[\lambda ]})\le (f: \phi )$.
\par Utilizing Formulas $(3.6.1)$, $(3.9.3)$
and $(3.9.8)$ we infer that
\par $(3.9.9)$ $(f_0: \phi )\le (f_0: f)(f^{[\lambda ]}: \phi ^{[\lambda ]})
\le (f_0: f)(f: \phi )$ for each $f_0\in \Upsilon (G,N_0)$ and
nonzero functions $f$ and $\phi $ in $C_{0,0}^+(G)$.
\par  Using $(3.5.1)$ and $(3.9.9)$ we infer that
$$(f_0: f)J_{\phi ,f_0}(f)=\frac{(f_0: f)(f:\phi )}{(f_0:\phi )
}\ge \frac{(f_0: \phi )}{(f_0:\phi ) }=1,$$ consequently, $J_{\phi
,f_0}(f)\ge (f_0: f)^{-1}$. Thus the left inequality in $(3.9.1)$
also is proved.
\par From Inequalities $(3.9.1)$ for $J_{\phi ,f_0}(f)$ and $J_{\phi ,f_0}(f_1)$
and Formula $(3.5.1)$ it follows $(3.9.2)$.

\par {\bf 3.10. Lemma.} {\it Let $G$ be a $T_1$ topological locally compact fan  loop
satisfying Condition $(3.5.2)$, let $f_0\in \Upsilon (G,N_0)$ (see
Condition $(3.5.3)$) and let $f_1$,...,$f_m$ be nonzero functions
belonging to $C_{0,0}^+(G)$, let also $0<\delta <\infty $, $0<\delta
_1<\infty $. Then there exists an open neighborhood $V$ of $e$ in
$G$ such that for each nonzero function $\phi $ in $C_{0,0}^+(G)$
with a support $S_{\phi }$ contained in $V$ and $0\le q_j\le \delta
_1$ for each $j=1,...,m$ the following inequality is satisfied:
$$(3.10.1)\quad \sum_{j=1}^mq_jJ_{\phi ,f_0}(f_j)\le J_{\phi
,f_0}(\sum_{j=1}^mq_jf_j)+\delta .$$}
\par {\bf Proof.} The loop $G$ is locally compact.
Let $S_{f_0,...,f_m}=\bigcup_{j=0}^m S_{f_j}$ be a common compact
support of these functions, where $S_{f_j}$ denotes a closed support
of $f_j$ (see also Subsection 3.1). We choose any function $g_1$ in
$C_{0,0}^+(G)$ such that $g_1: G\to [0,1]$ and
$g_1(S_{f_0,...,f_m}cl_G(Y_1))= \{ 1 \} $, where ${Y'}_1$ is an open
neighborhood of $e$ in $G$ with $Y_1=\check{Y'}_1$ and a compact
closure $cl_G(Y_1)$ (see Lemma 2.6). Consider arbitrary fixed
positive numbers $0<\delta <\infty $, $0<\delta _1<\infty $ and
$0<\epsilon <M$ such that \par $\epsilon \delta _1 \sum_{j=1}^m(f_j:
~ f_0 ) +\epsilon (1+\epsilon )(g_1: ~ f_0)\le \delta $, where
$M=\delta _1 m \max _{j=1,...m} \| f _j \|_G$. By virtue of
Corollary 2.15 the functions $f_0,...,f_m$ are uniformly $({\cal
L}_G, {\cal L}_H)$ continuous, where $H=({\bf C}, +)$. Therefore
there exists an open neighborhood $W'$ of $e$ with $W=\check{W'}$
and a compact closure $cl_G(W)$ in $G$ and $W\subset Y_1$, since $G$
is locally compact, such that
\par $(3.10.2)$ $|f_j(s)-f_j(x)|<\epsilon ^3[4Mm\delta _1]^{-1}$
\\ for each $s\setminus x\in W$.
Next we take a function $g\in C^+_{0,0}(G)$ such that $g: G\to
[0,1]$ and $g(S_{f_0,...,f_m}cl_G(W))= \{ 1 \} $ and $g(x)\le
g_1(x)$ for each $x\in G$, because $W\subset Y_1$. Hence $(g:
f_0)\le (g_1: f_0)$ by Inequality $(3.4.5)$.
\par  Let $S=P((S_{f_0,...,f_m}\cup
S_g)cl_G(W))$ (see Formula $(2.14.1)$). Since $cl_G(V)$,
$S_{f_0,...,f_m}$ and $S_g$ are compact, then $S$ is a compact
subset in $G$. For each open neighborhood $Y$ of $e$ in $G$ there
exists an open neighborhood $X$ of $e$ in $G$ such that $X^2\subset
Y$, since the multiplication in $G$ is continuous. In view of
Proposition 2.9 and Corollary 2.15 there exist open neighborhoods
${U'}_k$ of $e$ in $G$ such that $U_k=\check{U'}_k$ and
\par $(3.10.3)$ $[t((U_ka)U_k,(U_kb)U_k,(U_kc)U_k)U_k]
\cup [U_kt((U_ka)U_k,(U_kb)U_k,(U_kc)U_k)]$\par $\subset [t(a,b,c)
W_{k-1}]\cap [ W_{k-1}t(a,b,c)]$
\par $[p((U_ka)U_k,(U_kb)U_k,(U_kc)U_k)U_k]
\cup [U_kp((U_ka)U_k,(U_kb)U_k,(U_kc)U_k)]$\par $\subset [p(a,b,c)
W_{k-1}]\cap [ W_{k-1}p(a,b,c)]$
\\ for every $a$, $b$, $c$ in $S$ and $k\in \{ 1, 2 \} $ with
$U_0=W$ and an open of the form $(2.6.1)$ neighborhood $W_{k-1}$ of
$e$ in $G$ such that $\check{W}_{k-1}^2\subset U_{k-1}$ and
\par $(3.10.4)$ $|g(s)-g(x)|<\epsilon ^2[4M]^{-1}$ \\
for each $s$ and $x$ in $G$ such that $s\setminus x\in U_1$, where
$t=t_G$.
\par Take any $0\le q_j\le \delta _1$ for each
$j=1,...,m$ and put
\par $(3.10.5)$ $\Psi = \epsilon g + \sum_{j=1}^m q_jf_j$ and
\par $(3.10.6)$ $h_j(x)=q_jf_j(x)[\Psi (x)]^{-1}$ for each
$x\in S_{f_1,...,f_m}$ and \par $h_j(x)=0$ for each $x\in G-
S_{f_1,...,f_m}$, \\ where $S_{f_1,...,f_m}=\bigcup_{j=1}^m
S_{f_j}$. Therefore the function $\Psi $ belongs to $C^+_{0,0}(G)$
and $\sum_{j=1}^m h_j(x)\le 1$ for each $x\in G$.
\par From Inequalities $(3.10.2)$ and $(3.10.4)$ it follows that
\par $(3.10.7)$ $|\Psi (s)-\Psi (x)| \le \epsilon ^3[2M]^{-1}$
\\ for each $s$ and $x$ in $G$ such that $s\setminus x\in U_1$.
Moreover, $\| \Psi \|_G \le M+\epsilon <2M$.
\par Let $s$ and $x$ belong to $S_{f_1,...,f_m}cl_G(W)$
and $s\setminus x\in U_1$. The latter inclusion is equivalent to
$x\in sU_1$ and also to $s\in x/U_1$. Then from $(3.10.2)$ and
$(3.10.7)$ we deduce that
\par $(3.10.8)$ $|h_j(s)-h_j(x)|\le \epsilon /m$. \par Next we
consider the following case: $s\setminus x\in U_1$ and $x\notin
S_{f_1,...,f_m}cl_G(W)$. Suppose that $s\in S_{f_1,...,f_m}$, then
Condition $(3.10.3)$, Lemmas 2.2, 2.3 imply that $x\in
S_{f_1,...,f_m}cl_G(W)$ contradicting an assumption $x\notin
S_{f_1,...,f_m}cl_G(W)$. Hence $s\notin S_{f_1,...,f_m}$ and
consequently, $h_j(s)=0$ and $h_j(x)=0$. Thus Inequality $(3.10.8)$
takes place in this case as well. \par In the case $s\setminus x\in
U_1$ and $s\notin S_{f_1,...,f_m}cl_G(W)$ Condition $(3.10.3)$,
Lemmas 2.2, 2.3 imply that $x\notin S_{f_1,...,f_m}$. Therefore the
inequality $(3.10.8)$ is fulfilled in this case also. Thus the
estimate $(3.10.8)$ is satisfied for each $s$ and $x$ in $G$ such
that $s\setminus x\in U_1$.
\par Next we choose any fixed function $\phi \in C^+_{0,0}(G)$ such
that $\phi $ is not identically zero on $G$ and $\phi (y)=0$ for
each $y\in G- {U'}_2$. By virtue of Lemma 3.2 there are $m\in \bf
N$, $c_j>0$ and $b_j\in G$ for each $j\in \{ 1,...,m \} $ such that
$$(3.10.9) \quad \Psi (x)\le \sum_{j=1}^mc_j\phi (b_jx)$$ for each
$x\in G$ and
$$(3.10.10) \quad -\epsilon +\sum_{j=1}^mc_j\le (\Psi : ~ \phi )\le
\sum_{j=1}^mc_j.$$ Then Formulas $(3.10.3)$, $(3.10.8)$, $(3.10.9)$
and Lemma 2.2 imply that for each $x\in G$
$$\Psi (x)h_l(x)\le \sum_{j=1}^m c_j \phi (b_jx)[h_l(b_j\setminus e
) +\epsilon /m ]$$ for each $l$. Hence for each $x\in G$ we get
$$q_lf_l(x)=\Psi (x)h_l(x)\le \sum_{j=1}^m c_j
[h_l(b_j\setminus e ) +\epsilon /m ]\phi (b_jx)$$ and consequently,
$(q_lf_l: ~ \phi )\le \sum_{j=1}^m c_j [h_l(b_j\setminus e )
+\epsilon /m ]$. From $\sum_{l=1}^m h_l\le 1$ we deduce that
$\sum_{l=1}^m(q_lf_l: ~ \phi )\le (1+\epsilon ) \sum_{j=1}^m c_j$.
Together with Inequalities $(3.10.10)$ this leads to the following
estimate:
$$\sum_{j=1}^m (q_if_j: ~ \phi )\le (1+\epsilon ) (\Psi : ~ \phi ).$$
Dividing both sides of it on $(f_0: ~ \phi )$ we get the inequality
$$(3.10.11)\quad \sum_{j=1}^mq_jJ_{\phi ,f_0}(f_j)\le (1+\epsilon )
J_{\phi ,f_0}(\Psi ).$$ Then from $(3.7.1)$, $(3.7.2)$, $(3.10.5)$
and $(3.10.11)$ we infer that $$(3.10.12)\quad
\sum_{j=1}^mq_jJ_{\phi ,f_0}(f_j)\le J_{\phi
,f_0}(\sum_{j=1}^mq_jf_j)+\epsilon \sum_{j=1}^mq_jJ_{\phi ,f_0}(f_j)
+\epsilon (1+\epsilon )J_{\phi ,f_0}(g).$$ Therefore from
Inequalities $(3.9.1)$, $(3.10.12)$, $(3.4.5)$ and for $\epsilon $
as above it follows that
$$ \sum_{j=1}^mq_jJ_{\phi
,f_0}(f_j)\le J_{\phi ,f_0}(\sum_{j=1}^mq_jf_j)+\epsilon \delta _1
\sum_{j=1}^m(f_j: ~ f_0 )$$ $$ +\epsilon (1+\epsilon )(g: ~ f_0) \le
J_{\phi ,f_0}(\sum_{j=1}^mq_jf_j)+ \delta .$$ This implies the
estimate $(3.10.1)$ with $V={U'}_2$.

\par {\bf 3.11. Theorem.} {\it Let $G$ be a $T_1$ topological locally compact fan loop,
let $0<\epsilon $ and $f$ in $C_{0,0}^+(G)$ be a nonzero function,
$S_f=cl_G \{ x\in G: ~ f(x)\ne 0 \} $. Let also $V'$ be an open
neighborhood of $e$ in $G$ such that $V=\check{V'}$ and
\par $(3.11.1)$ $|f(x)-f(y)|<\epsilon $ for each $x$ and $y$ in $G$ with
$x\setminus y\in V$. Let $g\in C_{0,0}^+(G)$ be a nonzero function
such that $g(x)=0$ for each $x\in G- V'$. Then for each $\delta >
\epsilon $ and each open neighborhood ${W'}_e$ of $e$ in $G$ with
$W_e=\check{W'}_e$ and a compact closure $cl_G(W_e)$ contained in
$V$ there is an open neighborhood $U'$ of $e$ in $G$ such that
$U=\check{U'}$ and for each nonzero function $\phi $ in
$C_{0,0}^+(G)$ with a support $S_{\phi }$ contained in $U'$ there
are positive constants $c_1,...,c_n$ and elements $b_1,...,b_n$ in
$S_fcl_G(W_e)$ such that for each $x\in G$ and $\gamma \in N(G)$:
$$(3.11.2)\quad |f(\gamma x) - \sum_{j=1}^n \frac{c_j}{J^v_{\phi
,f_0}(g(v\setminus x)) } g(b_j\setminus \gamma x)|\le \delta ,$$
where an expression $J^v_{\phi ,f_0}(g(v\setminus x))$ means that a
functional $J_{\phi ,f_0}$ is taken in the $v$ variable.}
\par {\bf Proof.}  The continuous functions $f$
and $g$ are with compact supports, hence they are uniformly $({\cal
L}_G, {\cal L}_H)$ continuous and uniformly $({\cal R}_G, {\cal
R}_H)$ continuous on $G$ by Corollary 2.15, where $H=({\bf C}, +)$.
For each $y\in G$ the right translation operator $R_y$ is the
homeomorphism of $G$ as the topological space onto itself (see also
Subsection 2.5). Therefore a function $\nu (y) := (f(x): ~
g(x\setminus y))$ is continuous on the loop $G$ and consequently,
uniformly continuous on the compact subset $S_f$, hence $\sup_{y\in
S_f} \nu (y)< \infty $, where $(f(x): ~ g(x\setminus y))=(f: z)$ is
calculated in the $x$ variable with $z(x)=g(x\setminus y)$ for a
fixed parameter $y$. We take any fixed $\delta $ such that $\epsilon
<\delta <\infty $. Evidently there exists $0<\eta $ such that
\par $(3.11.3)$ $\eta \sup_{y\in S_f} \nu (y)<\delta - \epsilon $.
\par  Therefore take any fixed open
neighborhood ${W'}_e$ of $e$ in $G$ such that $W_e=\check{W'}_e$ and
$cl_G(W_e)$ is compact and $cl_G(W_e)\subset V$ (see Lemma 2.6). By
virtue of Corollary 2.15 the functions $g$ and $h$ are uniformly
$({\cal L}_G, {\cal L}_H)$ continuous and uniformly $({\cal R}_G,
{\cal R}_H)$ continuous. Hence there exists an open neighborhood
${W'}_1$ of $e$ in $G$ such that $W_1=\check{W'}_1$ and $cl_G(W_1)$
is compact and $cl_G(W_1)\subset {W'}_e$ and for each $x$ and $y$ in
$G$ with $x\setminus y\in W_1$:
\par $(3.11.4)$ $|g(x)-g(y)|<\eta $.
\par Therefore, a subset $S_f cl_G(W_1)$ is compact in $G$ (see Theorems
3.1.10, 8.3.13-8.3.15 in \cite{eng}, Lemma 2.6). Then we take
compact subsets $S_1=S_f cl_G(W_1)$ and $S=P(S_f cl_G(W_1))$ in $G$
(see Formula $(2.14.1)$). In view of Lemma 2.6 they contain open
subsets $S_fW_1$ and $P(S_fW_1)$ respectively, since $W_1$ is open
in $G$. Mention that the topological spaces $S_1$ and $S$ are
normal, since they are compact and $T_1\cap T_{3.5}$ (see Theorem
3.1.9 in \cite{eng}). Using Proposition 2.9 we take an open
neighborhood ${W'}_2$ of $e$ in $G$ with $W_2=\check{W'}_2$ such
that
\\ $(3.11.5)$ $[t((W_2a)W_2,(W_2b)W_2,(W_2c)W_2)W_2]\cup
[W_2t((W_2a)W_2,(W_2b)W_2,(W_2c)W_2)]$\par $ \subset
[t(a,b,c)W_3]\cap [W_3t(a,b,c)]$ and
\par $[p((W_2a)W_2,(W_2b)W_2,(W_2c)W_2)W_2]\cup
[W_2p((W_2a)W_2,(W_2b)W_2,(W_2c)W_2)] $\par $\subset
[p(a,b,c)W_3]\cap
[W_3p(a,b,c)]$ \\
for every $a$, $b$, $c$ in $S$, where $W_3$ is an open neighborhood
of $e$ in $G$ such that $\check{W}_3^2\subset W_1$.
\par In view of the Dieudonn\'e theorem 3.1 in \cite{hew} there
exists a partition of unity on $S_1$. Together with Theorem 3.3.2 in
\cite{eng} and Lemma 2.5 this implies that there are functions
$q_1,...,q_n$ in $C_{0,0}^+(G)$ and elements $w_1,...,w_n$ in $S_1$
such that $S_1\subset \bigcup_{j=1}^n w_jW_2$ and
\par $(3.11.6)$ $\sum_{j=1}^nq_j(x)=1$ for each $x\in S_1$ and
\par $(3.11.7)$ $q_j(y)=0$ for each $y\in G- (w_jW_2)$.
\par The conditions of this theorem imply that for each $x$
and $y$ in $G$ with $y\setminus x\in V$ the following inequalities
are satisfied: \par $(3.11.8)$ $[f(x)-\epsilon ]g(y\setminus x)\le
f(y)g(y\setminus x)\le [f(x) + \epsilon ]
g(y\setminus x)$, \\
since for $y\setminus x\in V$ Inequality $(3.11.1)$ is fulfilled;
for $u=y\setminus x\notin V$ the function $g$ is nil, $g(u)=0$.
\par Certainly $y\in w_jW_2$ if and only if there exists $b\in W_2$
such that $y=w_jb$. Then $(y\setminus x)\setminus (w_j\setminus
x)\in W_1$ if and only if there exists $c\in W_1$ such that
$w_j\setminus x=((w_jb)\setminus x)c$. For $w_j\setminus x=v\in V$
this gives $c=((w_jb)\setminus (w_jv))\setminus v$. In view of
$(2.2.2')$, $(2.2.4)$, $(2.1.8)$ and $(2.1.9)$ \par
$((w_jb)\setminus (w_jv))\setminus v=p(w_j,b,(w_jb)\setminus
(w_jv))((b\setminus v)\setminus v)$.
\par Therefore, from Conditions $(3.11.5)$-$(3.11.7)$ it follows that for each
$x$ and $y$ in $G$ and $j=1,...,n$:
\par $(3.11.9)$ $q_j(y)f(y)[g(y\setminus x)-\eta ] \le q_j(y)
f(y)g(w_j\setminus x) $\par $\le q_j(y)f(y)[g(y\setminus x)+\eta ]$. \\
Summing by $j$ in $(3.11.9)$, using $(3.11.8)$ we infer that for
each $x$ and $y$ in $G$: \par $(3.11.10)$ $[f(x)-\epsilon
]g(y\setminus x)- \eta f(y)$\par $\le \sum_{j=1}^n
q_j(y)f(y)g(w_j\setminus x)\le [f(x)+\epsilon ]g(y\setminus x)+ \eta
f(y)$.
\par Next we take any $\phi $ and $f_0$ in $C_{0,0}^+(G)$ such that
$\phi $ and $f_0$ are not identically zero. From Inequalities
$(3.11.10)$ after dividing on $J^y_{\phi ,f_0}(g(y\setminus x))$ and
Lemma 3.7 it follows that for each $x$ in $G$:
$$(3.11.11)\quad [f(x)-\epsilon ] -\eta \frac{ J_{\phi
,f_0}(f)}{J^y_{\phi ,f_0}(g(y\setminus x))}\le J^y_{\phi
,f_0}(\frac{\sum_{j=1}^ng(w_j\setminus x)q_j(y)f(y)}{J^v_{\phi
,f_0}(g(v\setminus x)) })$$ $$\le  [f(x)+\epsilon ] +\eta \frac{
J_{\phi ,f_0}(f)}{J^y_{\phi ,f_0}(g(y\setminus x))},$$ where
$J^y_{\phi ,f_0}(g(y\setminus u))=J_{\phi ,f_0}(z)$ means that the
functional $J^y_{\phi ,f_0}$ is taken in the $y$ variable in $G$,
where $z(y)=g(y\setminus u)$ for each $y\in G$ and a fixed parameter
$u$ in $G$.
\par Notice that the function $g(y\setminus x)$ is jointly
continuous in $(x,y)\in G\times G$. On the other hand, in view of
Lemmas 2.2, 2.4, 2.6 \par $\{ u=y\setminus x: ~ x\in S_f, ~ u\in S_g
\} $ is a compact subset in $G$, \\ since $Inv_l(S_f)$, $S_g$,
$S_fS_g$ and $t(S_f,Inv_l(S_f),S_fS_g)$ are compact subsets in $G$.
By virtue of Lemma 3.8 a mapping $\psi (x):=J^y_{\phi
,f_0}(g(y\setminus x))$ is continuous in the variable $x\in S_f$, $
~ \psi : S_f \to (0,\infty )$. Hence
\par $(3.11.12)$ $0<K_0=\inf_{x\in S_f}\psi (x) \le \sup_{x\in S_f}\psi (x)=K_1<\infty
$. \par Apparently in Formula $(3.11.3)$ the parameter $\eta >0$ can
be taken sufficiently small, because Inequalities $(3.11.3)$ and
$(3.11.12)$ are independent. Then from $(3.11.11)$ and $(3.11.12)$
we deduce that for each $\beta
>\epsilon $ there exist $q_j$ and $w_j$ (see above) such that \par
$\eta J_{\phi ,f_0}(f) <(\beta - \epsilon )\min (1, K_0)$,
consequently,
$$(3.11.13)\quad f(x)-\beta \le J^y_{\phi
,f_0}(\frac{\sum_{j=1}^ng(w_j\setminus x)q_j(y)f(y)}{J^v_{\phi
,f_0}(g(v\setminus x)) })\le  f(x)+\beta $$ for each $x\in G$.
\par In view of Lemmas 3.7 and 3.10 for each $\delta >\delta _1>\beta >\epsilon $
there exists an open of the form $(2.6.1)$ neighborhood $U$ of $e$
in $G$ such that $U\subset W_2$ and $$(3.11.14)\quad | J^y_{\phi
,f_0}(\frac{\sum_{j=1}^ng(w_j\setminus x)q_j(y)f(y)}{J^v_{\phi
,f_0}(g(v\setminus x)) }) -\sum_{j=1}^n \frac{J_{\phi
,f_0}(q_jf)}{J^v_{\phi ,f_0}(g(v\setminus x)) }g(w_j\setminus
x)|<\delta _1 - \beta
$$ for each $x\in S_f$. We put $c_j =J_{\phi ,f_0}(q_jf)$ and
$b_j=w_j$ for each $j=1,...,n$. Thus the estimates $(3.11.13)$ and
$(3.11.14)$ and Formula $(3.4.1')$ imply the assertion of this
theorem.

\par {\bf 3.12. Definition.} Let $W$ be an open neighborhood
of $e$ in a locally compact loop $G$ and a nonzero function $\phi
_{W}\in C^+_{0,0}(G)$ be such that $\phi _{W}(x)=0$ for each $x\in
G- W$. A family $ \{ \phi _W \} $ of these functions will be
directed by: \par $(3.12.1)$ $\phi _{W_1}\preceq \phi _{W_2}$ if and
only if $W_2\subseteq W_1$ and $\phi _{W_2}(x)=0$ implies $\phi
_{W_1}(x)=0$. If $\phi _{W_1}\preceq \phi _{W_2}$ and $\phi _{W_1}$
and $\phi _{W_2}$ are different functions, then it will be written
$\phi _{W_1}\prec \phi _{W_2}$.

\par {\bf 3.13. Lemma.} {\it Let $G$ be a $T_1$ topological locally compact fan  loop
satisfying Condition $(3.5.2)$ and let a family of nonzero functions
$ \{ \phi _U \} $ in $C^+_{0,0}(G)$ be directed by Condition
$(3.12.1)$. Let also $f_0\in \Upsilon (G,N_0)$ (see $(3.5.3)$) and
$f \in C^+_{0,0}(G)$. Then the limit exists: \par $(3.13.1)$
$\lim_{\{ \phi _U \} } J_{\phi _U, f_0}(f)=:J_{f_0}(f)$.}
\par {\bf Proof.} It is sufficient to prove that a net
$\{ J_{\phi _U, f_0}(f): \phi _U \} $ is fundamental (i.e. Cauchy)
in $\bf R$, where a net $ \{ \phi _U \} $ is directed by Condition
$(3.12.1)$. We take any fixed open neighborhood ${U'}_0$ of $e$ in
$G$ with $U_0=\check{U'}_0$ and a compact closure $cl_G(U_0)$. Let
$A=S_{f+f_0}cl_G(U_0)$, where $S_{f+f_0}=cl_G \{ x\in G: ~
f(x)+f_0(x)\ne 0 \} $. Therefore, a subset $S=P(A)$ is compact (see
Formula $(2.14.1)$ and Lemma 2.6), since $S_{f+f_0}$ is compact.
\par We choose any function $z\in C_{0,0}^+(G)$ such that $z|_A=1$.
Let $0<\epsilon <1$ and $\xi _1 = \epsilon (16[1+(z:f_0)]
[1+(f:f_0)])^{-1}$. From Corollary $(2.15)$ it follows that there
exists an open neighborhood $W'$ of $e$ in $G$ with $W=\check{W'}$
such that
\par $(3.13.2)$ $|f(x)-f(y)|<\xi _1 /2$ and
\par $(3.13.3)$ $|f_0(x)-f_0(y)|<\xi _1 /2$ \\
for each $x$ and $y$ in $G$ with $x\setminus y\in W$.
\par In view of Proposition 2.9 there exists an open
neighborhood ${U'}_2$ of $e$ in $G$ with $U_2=\check{U'}_2$ such
that
\par $(3.13.4)$ $[t((U_2a)U_2,(U_2b)U_2,(U_2c)U_2)U_2]\cup
[U_2t((U_2a)U_2,(U_2b)U_2,(U_2c)U_2)]$\par $ \subset
[t(a,b,c)B_1]\cap [B_1t(a,b,c)]$ and
\par $[p((U_2a)U_2,(U_2b)U_2,(U_2c)U_2)U_2]\cup
[U_2p((U_2a)U_2,(U_2b)U_2,(U_2c)U_2)]$\par $ \subset
[p(a,b,c)B_1]\cap [B_1p(a,b,c)]$ \\
for every $a$, $b$, $c$ in $S$, where $B_1$ is an open neighborhood
of $e$ in $G$ such that $\check{B}_1^2\subset U_1$, $U_1={U'}_0\cap
W'$ (see Lemma 2.6). Next we take a nonzero function $g\in
C_{0,0}^+(G)$ such that $g(x)=0$ for each $x\in G- {U'}_2$.
\par By virtue of Theorem 3.11
for any fixed $0<\delta <\xi _1 $ and each open neighborhood
${W'}_e$ of $e$ in $G$ with $W_e=\check{W'}_e$ and a compact closure
$cl_G(W_e)$ contained in ${U'}_2$ there is an open neighborhood
${U'}_{3,f}$ of $e$ in $G$ with $U_{3,f}=\check{U'}_{3,f}$ such that
for each nonzero function $\phi $ in $C_{0,0}^+(G)$ with a support
$S_{\phi }$ contained in ${U'}_{3,f}$ there are positive constants
$c_1,...,c_n$ and elements $b_1,...,b_n$ in $S_fcl_G(W_e)$ such that
for each $x\in G$ and $\gamma \in N(G)$:
$$(3.13.5)\quad |f(\gamma x)  - \sum_{j=1}^n \frac{c_j}{J^v_{\phi
,f_0}(g(v\setminus x)) } g(b_j\setminus \gamma x)|\le \delta .$$
Taking $U_{3,f}\subset {U'}_2$ we get $f(x)=0$ and $g(b_j\setminus
x)=0$ for each $x\in G- A$ according to the choice of $b_j$ in the
proof of Theorem 3.11, consequently,
$$(3.13.6)\quad |f(\gamma x) - \sum_{j=1}^n \frac{c_j}{J^v_{\phi
,f_0}(g(v\setminus x)) } g(b_j\setminus \gamma x)|\le \delta
z(\gamma x)$$ for each $x\in G$ and $\gamma \in N(G)$. From the
latter estimate and Lemma 3.7 we infer that
$$(3.13.7)\quad |J_{\phi ,f_0}(f) - K_{\phi ,f_0}(f;g)|\le
\delta J_{\phi ,f_0}(z)\le \delta (z: f_0), \mbox{  where}$$
$$K_{\phi ,f_0}(f;g)= J^x_{\phi ,f_0} (\sum_{j=1}^n
\frac{c_j}{J^v_{\phi ,f_0}(g(v\setminus x)) } g(b_j\setminus x)).$$
From Estimate $(3.13.7)$ and the right Inequality $(3.9.1)$ it
follows that
$$(3.13.7')\quad \sup_{\{ \phi _U \} }K_{\phi _U,f_0}(f;g)\le (1+\delta ) (f:f_0)+\delta
(z:f_0)<\infty .$$ Applying the proof above to $f_0$ instead of $f$
we get and open neighborhood ${U'}_{3,f_0}$ of $e$ with
$U_{3,f_0}=\check{U'}_{3,f_0}$ and $U_{3,f_0}\subset {U'}_2$ such
that for each nonzero function $\phi $ in $C_{0,0}^+(G)$ with a
support $S_{\phi }$ contained in ${U'}_{3,f_0}$  there are positive
constants $d_1,...,d_m$ and elements $v_1,...,v_m$ in
$S_{f_0}cl_G(W_e)$ such that
$$(3.13.8)\quad |f_0(\gamma x) - \sum_{j=1}^m \frac{d_j}{J^v_{\phi
,f_0}(g(v\setminus x)) } g(v_j\setminus \gamma x)|\le \delta
z(\gamma x)$$ for each $x\in G$ and $\gamma \in N(G)$, consequently,
$$(3.13.9)\quad |1- K_{\phi ,f_0}(f_0;g)|\le
\delta (z: f_0), \mbox{  where}$$
$$K_{\phi ,f_0}(f_0;g)= J^x_{\phi ,f_0} (\sum_{j=1}^m
\frac{d_j}{J^v_{\phi ,f_0}(g(v\setminus x)) } g(v_j\setminus x)),$$
since $J_{\phi ,f_0}(f_0)=1$. Moreover,
$$(3.13.10)\quad
\sup_{\{ \phi _U \} }K_{\phi _U,f_0}(f_0;g)\le (1+\delta )+ \delta
(z:f_0)<\infty .$$ Then ${U'}_3={U'}_{3,f}\cap {U'}_{3,f_0}$ is an
open neighborhood of $e$ in $G$. From $(3.13.7)$, $(3.13.9)$ and
$(3.13.10)$ we deduce that
$$(3.13.11)\quad |J_{\phi ,f_0}(f) -
\frac{K_{\phi ,f_0}(f;g)}{K_{\phi ,f_0}(f_0;g)}| \le \delta
_2+[1+\delta +\delta _2]\delta _2(1-\delta _2)^{-1},$$ where $\delta
_2=\delta (z: f_0)<\xi _1 (z: f_0)<1/16$.  In view of Lemmas 3.7 and
3.10, Formulas $(3.13.5)$ and $(3.13.6)$ there exists an open
neighborhood ${U'}_4$ of $e$ with $U_4=\check{U'}_4$ and $U_4$
contained in ${U'}_3$ such that for each nonzero $\phi $ in
$C_{0,0}^+(G)$ with $S_{\phi }\subset {U'}_4$ there are
inequalities:
$$(3.13.12)\quad |K_{\phi ,f_0}(f;g)- \sum_{j=1}^nc_jJ^x_{\phi ,f_0}
(\frac{g(b_j\setminus x)}{J^v_{\phi ,f_0}(g(v\setminus \gamma x))})
| \le \delta \mbox{ and}$$
$$(3.13.13)\quad |K_{\phi ,f_0}(f_0;g)-
\sum_{j=1}^md_jJ^x_{\phi ,f_0}(\frac{g(v_j\setminus x)}{J^v_{\phi
,f_0}(g(v\setminus \gamma x))}) | \le \delta $$ for each $\gamma \in
N(G)$. On the other hand, Formulas $(3.5.1)$, $(3.6.4)$, $(3.6.4')$,
$(3.6.12)$ and $(3.4.3)$ imply that
$$(3.13.14)\quad J^x_{\phi ,f_0}(\frac{g(b_j\setminus x)}{J^v_{\phi ,f_0}(g(v\setminus \gamma x))})=
\int_{N_0}(\frac{g(b_j\setminus x)}{(g(v\setminus \gamma x):\phi
(v))}: \phi (x))\lambda (d\gamma )$$ $$=
 \frac{(g(b_j\setminus x): \phi (x))}{(g^{[\lambda ]}(v\setminus e): \phi
(v))}.$$ Then from Proposition $(2.9.1)$ and Formulas $(3.4.2)$,
$(3.4.2')$ it follows that for each $b\in G$ and each $0<\delta
_3\le \delta $ there exists an open neighborhood ${U'}_{5,b}$ of $e$
in $G$ with $U_{5,b}=\check{U'}_{5,b}$ such that for each nonzero
$\phi _U\in C_{0,0}^+(G)$ with $S_{\phi _U}\subset U\subset
{U'}_{5,b}$
$$(3.13.15)\quad |\frac{(g(b\setminus x):\phi
_U(x))}{(g^{[\lambda ]}(v\setminus e):\phi _U(v))}-\frac{(g(x): \phi
_U(x))}{(g^{[\lambda ]}(v\setminus e): \phi _U(v))}|<\delta _3,$$
since $S_{\phi _U}\subset U$ and $t(a,b,e)=t(a,e,b)=t(e,a,b)=e$ and
$p(a,b,e)=p(a,e,b)=p(e,a,b)=e$ for each $a$ and $b$ in $G$.
Therefore we take ${U'}_5=\bigcap_{j=1}^n {U'}_{5,b_j}\cap
\bigcap_{k=1}^m{U'}_{5,v_k}\cap {U'}_4$ and $\phi =\phi _Y$ with
$Y={U'}_5$. We put $c=\sum_{j=1}^nc_j$ and $d=\sum_{k=1}^md_k$. From
$(3.13.11)-(3.13.15)$ and $(3.9.1)$ it follows that
$$\frac{c}{d}< K_1, \mbox{  where   }K_1=3[1+(f:f_0)] (1+\delta
)(1-\delta )^{-1}<4[1+(f:f_0)].$$  Then we deduce from Formulas
$(3.13.11)$-$(3.13.15)$ for each $\phi _U$ with an open neighborhood
$U$ of $e$ in $G$ such that $U\subset {U'}_5$:
$$|J_{\phi _U,f_0}(f) - \frac{c}{d}|< \delta (1-\delta
)^{-1}4[1+(f:f_0)]+\delta _2+[1+\delta +\delta _2]\delta _2
(1-\delta _2)^{-1},$$ consequently, \par $(3.13.16)$ $|J_{\phi
_{V_1},f_0}(f) - J_{\phi _{V_2},f_0}(f) |<8 \delta (1-\delta
)^{-1}[1+(f:f_0)]$\par $+2\delta _2+2[1+\delta +\delta _2]\delta _2(1-\delta _2)^{-1}<\epsilon $ \\
for each open neighborhoods $V_1$ and $V_2$ of $e$ in $G$ such that
$V_1 \subset {U'}_5$ and $V_2\subset {U'}_5$. Thus the net $ \{
J_{\phi _U,f_0}(f) : \phi _U \} $ is fundamental, where the net $ \{
\phi _U \} $ is directed by Condition $(3.12.1)$.

\par {\bf 3.14. Remark.} Suppose that $G$ is a $T_1$
topological locally compact fan  loop and Condition $(3.5.2)$ is
fulfilled and $f_0\in \Upsilon (G,N_0)$ (see $(3.5.3)$), functions
$f$ and $g$ belong to $C^+_{0,0}(G)$ and $g$ is nonzero. Then in
view of Lemma 3.13 a functional exists \par $(3.14.1) \quad
J_{g}(f)=J_{f_0}(f)/J_{f_0}(g)$.
\par As a consequence of Lemma 3.13 and
Formulas $(3.5.1)$ and $(3.14.1)$ we get that
\par $(3.14.2)$ the functional $J_{g}(f)$ is independent of $f_0$.
\par Then Formula $(3.9.2)$ and Lemma 3.13 imply that
\par $(3.14.3) \quad  (g: f_0)^{-1}(f_0: f )
^{-1}\le J_{g}(f) \le (f: f_0) (f_0: g)$ \\ for each $f_0 \in
\Upsilon (G,N_0)$ and a nonzero function $f\in C^+_{0,0}(G)$.

\par {\bf 3.15. Theorem.} {\it Let $G$ be a $T_1$
topological locally compact fan  loop fulfilling Condition $(3.5.2)$
and a functional $J=J_g$ be defined by Formula $(3.14.1)$. Then $J$
possesses the following properties:
\par $(3.15.1)$ $J(f)\ge 0$ for each $f\in C^+_{0,0}(G)$; and
 if a function $f$ is nonzero, then $J(f)>0$;
 \par $(3.15.2)$ $J(\alpha _1f_1+...+\alpha _nf_n)=
 \alpha _1J(f_1)+...+\alpha _nJ(f_n)$
 for each $f_1,...,f_n$ in $C^+_{0,0}(G)$ and $\alpha _1\ge
 0$,...,$\alpha _n\ge 0$;
\par $(3.15.3)$ $J(\mbox{}_bf)=J(f)$ for each $b\in G$ and $f\in C^+_{0,0}(G)$.}
\par {\bf Proof.} Property $(3.15.1)$ follows
from Formula $(3.14.3)$. On the other hand, Lemmas 3.7, 3.10, 3.13
imply Equality $(3.15.2)$.
\par Then Formulas $(3.6.4)$, $(3.6.4')$, $(3.6.12)$, $(3.14.1)$ and Lemma 3.13
imply
\par $(3.15.4)$ $J(\mbox{}_bf^{[\lambda ]})=J(f^{[\lambda ]})$ for each $b\in G$
and $f$ in $C^+_{0,0}(G)$. \par As a topological space $G$ is
locally compact. According to the measure theory on locally compact
spaces (see Chapter 3, Section 11 in \cite{hew}) a functional $J$ on
$C^+_{0,0}(G)$ satisfying Conditions $(3.15.1)$ and $(3.15.2)$
induces
\par $(3.15.5)$ a regular $\sigma $-additive measure $\mu $ on
a Borel $\sigma $-algebra ${\cal B}(G)$ of $G$ such that $\mu (U) =
\sup \{ \mu (X): ~ X \mbox{ is compact}, ~ X\subset U \} $ \\ for
each open subset $U$ in $G$ and
\par $(3.15.6)$ $\mu (A) = \inf \{ \mu (V): ~ V \mbox{ is open},
 ~ A\subset V\subset G \} $ \\ for each $A\in {\cal B}(G)$ and
\par $(3.15.7)$ $J(f)=\int_G f(x) \mu (dx)$ for each $f\in
C^+_{0,0}(G)$ and \par $(3.15.8)$ the functional $J$ has an
extension $\bar{\bar{J}}$ such that $\bar{\bar{J}}(f)=\int_Gf(x)\mu
(dx)$ for each nonnegative $\mu $-measurable function $f$ on $G$,
where\par $\bar{\bar{J}}(f)=\inf \{ \bar{J}(h): ~ h\ge f, ~ h \mbox{
is lower semicontinuous } \} ,$ \par $\bar{J}(h)= \sup \{ J(p): ~
p\in C^+_{0,0}(G), ~ p\le h \} $ \\ (see Theorems 11.22, 11.23,
11.36 and Corollary 11.37 in \cite{hew}).
\par On the other hand, for each $\gamma \in N(G)$ Formulas $(3.4.1')$ and $(3.4.2')$ give
\par $(3.15.9)$ $(\mbox{}_{\gamma }f: \phi _U)= (f: \mbox{}_{\gamma }\phi _U)=(f:
\phi _U)$. \par  From Lemma 3.13, Formulas $(3.14.1)$ and $(3.15.9)$
we deduce that
\par $(3.15.10)$ $J(\mbox{}_{\gamma }f)=J(f)$ for each
$\gamma \in N_0(G)$.
\par By virtue of the Fubini theorem 13.8 in \cite{hew}, $(3.5.3)$, $(3.5.4)$,
Formulas $(3.15.4)$, $(3.15.7)$ and $(3.15.10)$ above we infer that
$$J(\mbox{}_bf)=\int_{N_0} J(\mbox{}_{b\gamma }f) \lambda (d\gamma
)= \int_G \int_{N_0} \mbox{ }_bf(\gamma x)\lambda (d\gamma )\mu
(dx)$$ $= J(\mbox{}_bf^{[\lambda ]})=J(f^{[\lambda
]})=\int_{N_0}J(\mbox{}_{\gamma }f)\lambda (d\gamma )=J(f)$,
\\ since $\lambda (N_0)=1$ and $N_0\subset N(G)$. Thus the last assertion of this theorem
also is proved.

\par {\bf 3.16. Theorem.} {\it If $G$ is a
$T_1$ topological locally compact fan  loop fulfilling Condition
$(3.5.2)$, then there exists
\par $(3.16.1)$ a regular $\sigma $-additive measure $\mu $ on a
Borel $\sigma $-algebra ${\cal B}(G)$ of $G$, $\mu : {\cal B}(G)\to
[0, \infty ]$ such that
\par $(3.16.2)$ $\mu (U)>0$ for each open subset in $G$;
\par $(3.16.3)$ $\mu (A)<\infty $ for each compact subset $A$ in
$G$; \par $(3.16.4)$ $\mu (bB)=\mu (B)$ for each $B\in {\cal B}(G)$
and $b\in G$. \par Such $\mu $ can be chosen corresponding to a
functional $J$ satisfying Conditions $(3.15.1)$-$(3.15.3)$.}
\par {\bf Proof.} This is an immediate consequence of
$(3.15.1)$-$(3.15.3)$, $(3.15.5)$-$(3.15.8)$. In particular $\mu
(A)=\bar{\bar{J}}(\chi _A)$ for the characteristic function $\chi
_A$ of a Borel subset $A$ in $G$, where $\chi _A(x)=1$ for each
$x\in A$, $\chi _A(y)=0$ for each $y\in G-A$.

\par {\bf 3.17. Remark.} Each function $f$ in $C_{0,0}(G)$ can be
represented as $f=f^+ - f^-$, where $f^+(x)=\max (0, f(x))$, $f^+$
and $f^-$ belong to $C^+_{0,0}(G)$. Therefore, a functional $J$
satisfying Conditions $(15.1)$ and $(15.2)$ can be extended to a
linear functional on $C_{0,0}(G)$ such that $J(f)=J(f^+)-J(f^-)$.
Hence Property $(3.15.3)$ extends onto $C_{0,0}(G)$.

\par {\bf 3.18. Definition.} A linear functional $J$ on
$C_{0,0}(G)$ satisfying Property $(3.15.3)$ is called left
invariant.
\par A measure $\mu $ on the Borel $\sigma $-algebra ${\cal B}(G)$ of a
topological fan  loop $G$ such that $\mu $ satisfies Condition
$(3.16.4)$ is called left invariant.

\par {\bf 3.19. Theorem.} {\it Let $G$ be a
$T_1$ topological locally compact fan  loop fulfilling Condition
$(3.5.2)$ and let $\mu $ be a measure possessing Properties
$(3.16.1)$-$(3.16.4)$. Then $\mu (G)<\infty $ if and only if $G$ is
compact.}
\par {\bf Proof.} If $G$ is compact, then by $(3.16.3)$
$\mu (G)<\infty $.
\par Vice versa suppose that $\mu (G)<\infty $ and
consider the variant that $G$ is not compact and take an open
neighborhood $U'$ of $e$ in $G$ with $U=\check{U}'$ such that
$U=N_0U$ and its closure $cl_G(U)$ is compact, hence $0<\mu
(U)<\infty $ (see also Condition $(3.5.2)$). By virtue of Theorem
2.8 there exists an open neighborhood $V'$ of $e$ in $G$ with
$V=\check{V}'$ such that $V=N_0V$ and $[cl_G(V)]^2\subset U'$. In
view of Lemma 2.5 a subset $xU$ is open in $G$ for each $x\in G$.
\par At first we take some fixed $x_1\in G$. Then we construct a
sequence $\{ x_j: ~ j \in {\bf N} \} $ by induction. Let
$x_1,...,x_n$ be constructed such that if $n\ge 2$, then $x_jV\cap
x_kV=\emptyset $ for each $1\le j<k\le n$. There exists \par $y\in
G-\bigcup_{j=1}^nU_j$, where
$U_j:=x_jUp(x_jU,V,V)p(V,V,V)[p(x_jU,V,V)]^{-1}$,
\\ since $G$ is not compact and $U_j$ is open by Lemma 2.6 and $cl_G(U_j)$
is compact by Theorem 3.1.10 in \cite{eng} and Lemmas 2.4, 2.6. Put
$x_{n+1}=y$ with this $y$. \par Suppose that there is $z\in x_jV\cap
x_{n+1}V$ for some $1\le j\le n$. Therefore there would be $v$ and
$u$ in $V$ for which $z=x_jv=x_{n+1}u$, consequently,
$(x_jv)/u=(x_{n+1}u)/u=x_{n+1}$ by Condition $(2.1.2)$ and Formula
$(2.2.5)$. Therefore by Formulas $(2.2.3)$, $(2.3.3')$ and Condition
$(2.1.9)$
\par $x_{n+1}=x_j(v(e/u))p(x_j,v,e/u)p(e/u,u,u\setminus e)
[p(x_j(v(e/u)),u,u\setminus e)]^{-1}$\\ contradicting the choice of
$x_{n+1}$, since $[cl_G(V)]^2\subset U'$. Thus $x_jV\cap
x_kV=\emptyset $ for each $1\le j<k\le n+1$. This would mean by
$(3.16.4)$ that $\mu (G)\ge \sum_{j=1}^n \mu (x_jV)=n\mu (V)$ for
each $n$, contradicting $0<\mu (G)<\infty $.

\par {\bf 3.20. Theorem.} {\it Assume that $G$ is a
$T_1$ topological locally compact fan  loop satisfying Condition
$(3.5.2)$ and functionals $J$ and $H$ on $C^+_{0,0}(G)$ satisfy
Conditions $(3.15.1)$-$(3.15.3)$.
\par Then a positive constant $\kappa $ exists such that
\par $(3.20.1)$ $H(f)=\kappa J(f)$ for each $f\in C^+_{0,0}(G)$.}
\par {\bf Proof.} By virtue of Theorem 3.16 there exist two measures
$\mu _1$ and $\mu _2$ corresponding to $J$ and $H$. We consider a
subalgebra ${\cal C}(G):=\theta ^{-1}({\cal B}(G/\cdot /N_0))$ in
${\cal B}(G)$, where $\theta : G\to G/\cdot /N_0$ is the quotient
homomorphism, ${\cal B}(G)$ denotes the Borel $\sigma $-algebra on
$G$. Put $\nu _j(A)= \mu _j(\theta ^{-1}(A))$ for each $j$ and $A\in
{\cal B}(G/\cdot /N_0)$. \par From Theorems 2.8 and 3.16 it follows
that the measure $\nu _j$ on the group $G/\cdot /N_0$ is such that
$\nu _j(V)>0$ for each open subset $V$ in $G/\cdot /N_0$, $\nu
_j(A)<\infty $ for each compact subset $A$ in $G/\cdot /N_0$, $\nu
_j(cB)=\nu _j(B)$ for each $c\in G/\cdot /N_0$ and $B\in {\cal
B}(G/\cdot /N_0)$, $j \in \{ 1, 2 \} $. By virtue of Theorem 15.6 in
\cite{hew} there are positive constants $p_j$ such that $\nu
_j=p_j\eta $, where $\eta $ is a left invariant Haar measure on
$G/\cdot /N_0$. Thus $J(f^{[\lambda ]})=p_1H(f^{[\lambda ]})/p_2$
for each $f\in C^+_{0,0}(G)$.
\par We consider $\eta _1(b,f)=J(\mbox{}_bf)/J(f^{[\lambda ]})$ and
$\eta _2(b,f)=H(\mbox{}_bf)/H(f^{[\lambda ]})$ for each $b\in G$ and
a nonzero function $f$ in $C^+_{0,0}(G)$. According to Property
$(3.15.3)$ we get the identities $\eta _j(b,f)=\eta _j(e,f^{[\lambda
]})=1$ for each $j\in \{ 1, 2 \} $. This implies that for each
nonzero function $f\in C^+_{0,0}(G)$ and $b\in G$:
\par $(3.20.2)$ $J(\mbox{}_bf)/H(\mbox{}_bf)=p_1/p_2$. \par The
measures $\mu _1$ and $\mu _2$ possess Properties
$(3.16.1)$-$(3.16.4)$. In view of the Lebesgue-Radon-Nikodym theorem
(see Theorem (12.17) in \cite{hew} or see \cite{bogachmtb}) there
exists a $\mu _1$ measurable nonnegative function $h(x)$ such that
$\int_Gg(x)\mu _2(dx)=\int_Gg(x)h(x)\mu _1(dx)$ for each $g\in
C^+_{0,0}(G)$. Therefore from Formulas $(3.15.8)$ and $(3.20.2)$ it
follows that $h(x)$ is a positive constant. Thus $(3.20.1)$ is
proved.

\section{Appendix. Products of fan loops.}
\par The main subject of this paper are measures
on fan loops. Nevertheless, in this section it is shortly
demonstrated that there are abundant families of fan loops besides
those which appear in areas described in the introduction.

\par {\bf 4.1. Theorem.} {\it Let $(G_j,\tau _j)$ be a family of
topological $T_1$ fan loops (see Definition 2.1), where $j\in J$,
$J$ is a set. Then their direct product $G=\prod_{j\in J}G_j$
relative to the Tychonoff product topology $\tau $ is a topological
$T_1$ fan loop and
\par $(4.1.1)$ $Z(G)=\prod_{j\in J}Z(G_j)$ and $N(G)=\prod_{j\in J}N(G_j)$.}
\par {\bf Proof.} The direct product of topological loops
is a topological loop (see \cite{bruckb,eng,kakkar}). Thus
conditions $(2.1.1)$-$(2.1.3)$ are satisfied. \par Each element
$a\in G$ is written as $a= \{ a_j: ~ \forall j\in J, ~ a_j\in G_j
\}$. From $(2.1.4)$-$(2.1.7)$ we infer that \par $(4.1.2)$ $Com (G)
:= \{ a\in G: \forall b\in G, ~ ab=ba \} =$\par $ \{ a\in G: ~ a= \{
a_j: \forall j\in J, a_j\in G_j \}; \forall b\in G, ~ b= \{ b_j:
\forall j\in J, b_j\in G_j \} ; \forall j\in J, ~ a_jb_j=b_ja_j \}
=\prod_{j\in J} Com (G_j)$,
\par $(4.1.3)$ $N_l(G) := \{a\in G: ~ \forall b\in G, ~ \forall c\in G, ~ (ab)c=a(bc)
\} = \{a\in G: ~ a= \{ a_j: \forall j\in J, a_j\in G_j \}; ~ \forall
b\in G, ~ b= \{ b_j: \forall j\in J, b_j\in G_j \}; ~ \forall c\in
G, ~ c= \{ c_j: \forall j\in J, c_j\in G_j \}; ~ \forall j\in J,~
(a_jb_j)c_j=a_j(b_jc_j) \}= \prod_{j\in J} N_l(G_j)$
\\ and similarly \par $(4.1.4)$ $N_m(G)=\prod_{j\in J} N_m(G_j)$ and
\par $(4.1.5)$ $N_r(G)=\prod_{j\in J} N_r(G_j)$. \\ Therefore $(4.1.3)$-$(4.1.5)$
and $(2.1.8)$ imply that
\par $(4.1.6)$ $N(G)=\prod_{j\in J}N(G_j)$. Thus
\par $(4.1.7)$ $Z(G) := Com (G)\cap N(G)=\prod_{j\in J}Z(G_j)$.
\par Let $a$, $b$ and $c$ be in $G$, then
\par $(ab)c=\{ (a_jb_j)c_j: ~ \forall j \in J, ~ a_j\in G_j, b_j\in
G_j, c_j\in G_j \} $\par $= \{ t_{G_j}(a_j,b_j,c_j) a_j(b_jc_j): ~
\forall j \in J, ~ a_j\in G_j, b_j\in G_j, c_j\in G_j \}$\par $ =
t_G(a,b,c) a(bc)$ \par and analogously $(ab)c = a(bc)p_G(a,b,c) $,
where
\par $(4.1.8)$ $t_G(a,b,c) = \{ t_{G_j}(a_j,b_j,c_j): ~ \forall j \in J, ~
a_j\in G_j, b_j\in G_j, c_j\in G_j \} $ and
\par $(4.1.9)$ $p_G(a,b,c) = \{ p_{G_j}(a_j,b_j,c_j): ~ \forall j \in J, ~
a_j\in G_j, b_j\in G_j, c_j\in G_j \} $.
\par Therefore, Formulas
$(4.1.7)$-$(4.1.9)$ imply that Conditions $(2.1.9)$ also are
satisfied. Thus $G$ is a topological fan loop. By virtue of Theorem
2.3.11 in \cite{eng} a product of $T_1$ spaces is a $T_1$ space,
hence $G$ is the $T_1$ topological fan loop.

\par {\bf 4.2. Corollary.} {\it $(1)$. Let conditions of Theorem 4.1 be satisfied
and for each $j\in J$ a fan loop $G_j$ satisfies Condition
$(3.5.2)$. Then the product fan loop $G$ satisfies Condition
$(3.5.2)$. \par $(2)$. Moreover, if $G_j$ is compact for all $j\in
J_0$ and locally compact for each $j\in J\setminus J_0$, where
$J_0\subset J$ and $J\setminus J_0$ is a finite set, then $G$ is
locally compact.}
\par {\bf Proof.} Using Formulas $(4.1.8)$ and $(4.1.9)$ it is sufficient to take
$N_0(G)=\prod_{j\in J} N_0(G_j)$, since the direct product of
compact groups $N_0(G_j)$ is a compact group $N_0(G)$ (see the
Tychonoff theorem 3.2.4 in \cite{eng} or \cite{hew}). The last
assertion $(2)$ follows from the known fact that $G$ as a
topological space is locally compact under the imposed above
conditions (see Theorem 3.3.13 in \cite{eng}).

\par {\bf 4.3. Remark.}
\par $(4.3.1)$. Let $A$ and $B$ be two fan loops
and let $N$ be a group such that $N_0(A)\hookrightarrow N$,
$N_0(B)\hookrightarrow N$, $N\hookrightarrow N(A)$ and
$N\hookrightarrow N(B)$ and let $N$ be normal in $A$ and in $B$ (see
also Sections 2.1, 2.7 and 3.5).
\par Using direct products it is always possible to extend either
$A$ or $B$ to get such a case. In particular, either $A$ or $B$ may
be a group. On $A\times B$ an equivalence relation $\Xi $ is
considered such that
\par $(4.3.2)$ $(v\gamma ,b)\Xi (v,\gamma b)$
\\ for every $v$ in $A$, $b$ in $B$ and $\gamma $ in $N$.
\par $(4.3.3)$. Let $\phi : A\to {\cal A}(B)$ be a
single-valued mapping, where ${\cal A}(B)$ denotes a family of all
bijective surjective single-valued mappings of $B$ onto $B$
subjected to the conditions given below. If $a\in A$ and $b\in B$,
then it will be written shortly $b^a$ instead of $\phi (a)b$, where
$\phi (a) : B\to B$. Let also \par $\eta _{\phi }: A\times A\times
B\to N$, $\kappa _{\phi } : A\times B\times B\to N$ \par and $\xi
_{\phi } : ((A\times B)/\Xi ) \times ((A\times B)/\Xi )\to N$ \\
be single-valued mappings written shortly as $\eta $, $\kappa $ and
$\xi $ correspondingly such that
\par $(4.3.4)$ $(b^u)^v=b^{vu}\eta (v,u,b)$, $~ {\gamma }^u=\gamma $, $~b^{\gamma }=b$;
\par $(4.3.5)$ $\eta (v,u,(\gamma _1b)\gamma _2)=\eta (v,u,b)$;
\par if $\gamma \in \{ v, u, b \} $ then $\eta (v,u,b)=e$;
\par $(4.3.6)$ $(cb)^u=c^ub^u\kappa (u,c,b)$;
\par $(4.3.7)$ $\kappa (u,(\gamma _1c)\gamma _2,(\gamma _3b)\gamma _4)=\kappa
(u,c,b)$ and \par if $\gamma \in \{ u, c, b)$ then $\kappa (u,c
,b)=e$;
\par $(4.3.8)$ $\xi (((\gamma u)\gamma _1,(\gamma _2c)\gamma _3),((\gamma _4v)\gamma _5,
(\gamma _6b)\gamma _7))= \xi ((u,c),(v,b))$ and
\par $\xi ((e,e), (v,b))=e$ and $\xi ((u,c),(e,e))=e$
\\ for every $u$ and $v$ in $A$, $b$, $c$ in $B$,
$\gamma $, $\gamma _1$,...,$\gamma _7$ in $N$, where $e$ denotes the
neutral element in $N$ and in $A$ and $B$.
\par We put
\par $(4.3.9)$ $(a_1,b_1)(a_2,b_2)=(a_1a_2,b_1b_2^{a_1}\xi
((a_1,b_1),(a_2,b_2)))$ \\ for each $a_1$, $a_2$ in $A$, $b_1$ and
$b_2$ in $B$. \par The Cartesian product $A\times B$ supplied with
such a binary operation $(4.3.9)$ will be denoted by $A\bigotimes
^{\phi , \eta , \kappa , \xi }B$.

\par {\bf 4.4. Theorem.} {\it Let the conditions of Remark 4.3 be
fulfilled. Then the Cartesian product $A\times B$ supplied with a
binary operation $(4.3.9)$ is a fan loop.}
\par {\bf Proof.} From the conditions of Remark 4.3 it follows
that the binary operation $(4.3.9)$ is single-valued. The group $N$
is normal in the loops $A$ and $B$ by Conditions $(4.3.1)$. Hence
for each $a\in A$ and $\beta \in N$ there exists $(a\beta )/a \in N$
and $a\setminus (\beta a)\in N$, since $aN=Na$ for each $a\in A$.
Similarly it is for $B$. Thus there are single-valued mappings \par
$r_{A,a}(\beta )=(a\beta )/a$, $ ~ \check{r}_{A,a}(\beta
)=a\setminus (\beta a)$,
\par $r_{B,b}(\beta )=(b\beta )/b$, $ ~ \check{r}_{B,b}(\beta
)=b\setminus (\beta b)$, \par $r_{A,a}: N\to N$, $ ~
\check{r}_{A,a}: N\to N$, $~ r_{B,b}: N\to N$, $ ~ \check{r}_{B,b}:
N\to N$ \\ for each $a\in A$ and $b\in B$. Evidently \par $r_{A,a}(
\check{r}_{A,a}(\beta ))=\beta $ and $\check{r}_{A,a}(r_{A,a}(\beta
))=\beta $ \\ for each $a\in A$ and $\beta \in N$, and similarly for
$B$.
\par Let $I_1=((a_1,b_1)(a_2,b_2))(a_3,b_3)$ and $I_2=
(a_1,b_1)((a_2,b_2)(a_3,b_3))$, where $a_1$, $a_2$, $a_3$ belong to
$A$, $b_1$, $b_2$, $b_3$ belong to $B$. Then we infer that
\par $I_1= ((a_1a_2)a_3,(b_1b_2^{a_1})\xi ((a_1,b_1),(a_2,b_2))
b_3^{a_1a_2}\xi ((a_1a_2,b_1b_2^{a_1}),(a_3,b_3)))$ and
\par $I_2= (a_1(a_2a_3), b_1(b_2^{a_1}b_3^{a_1a_2})\beta )$ with \par $\beta =\eta
(a_1,a_2,b_3)\kappa (a_1,b_2,b_3^{a_2})[\xi
((a_2,b_2),(a_3,b_3))]^{a_1} \xi ((a_1,b_1),
(a_2a_3,b_2b_3^{a_2}))$. Hence \par $I_1=(a,b\alpha )$ and
$I_2=(a,b\beta )$, where  $a=a_1(a_2a_3)$ and
$b=b_1(b_2^{a_1}b_3^{a_1a_2})$,
\par $\alpha =\check{r}_{B,b}(
p_A(a_1,a_2,a_3))$\par $p_B(b_1,b_2^{a_1},b_3^{a_1a_2}) ~
\check{r}_{B,b_3^{a_1a_2}}(\xi ((a_1,b_1),(a_2,b_2))) ~ \xi
((a_1a_2,b_1b_2^{a_1}),(a_3,b_3)))$. \par Therefore
\par $(4.4.1)$ $I_1=I_2p$ with $p=p_{A\bigotimes ^{\phi , \eta , \kappa , \xi }B}((a_1,b_1),(a_2,b_2),(a_3,b_3))$
and  \par $I_1=tI_2$ with $t=t_{A\bigotimes ^{\phi , \eta , \kappa ,
\xi }B}((a_1,b_1),(a_2,b_2),(a_3,b_3))$;
\par $(4.4.2)$ $p=\beta ^{-1}\alpha $ and $t=r_{A,a}(r_{B,b}(p))$.
\\ Apparently
$t_{A\bigotimes ^{\phi , \eta , \kappa , \xi }B}
((a_1,b_1),(a_2,b_2),(a_3,b_3))\in N$ and \\ $p_{A\bigotimes ^{\phi
, \eta , \kappa , \xi }B} ((a_1,b_1),(a_2,b_2),(a_3,b_3))\in N$ for
each $a_j\in A$, $b_j\in B$, $j\in \{ 1, 2, 3 \} $, since $\alpha $
and $\beta $ belong to the group $N$.

\par If $\gamma \in N$ and either $(\gamma ,e)$ or $(e,\gamma )$ belongs
to $\{ (a_1,b_1),(a_2,b_2),(a_3,b_3) \} $, then from the conditions
of Section 4.3 and Formulas $(4.4.1)$ and $(4.4.2)$ it follows that
\par $p_{A\bigotimes ^{\phi , \eta , \kappa , \xi }B}((a_1,b_1),(a_2,b_2),(a_3,b_3))=e$ and
\par  $t_{A\bigotimes ^{\phi , \eta , \kappa , \xi
}B}((a_1,b_1),(a_2,b_2),(a_3,b_3))=e$, \\ consequently, $(N,e)\cup
(e,N)\subset N(A\bigotimes ^{\phi , \eta , \kappa , \xi }B)$.

\par Apparently $(2.1.3)$ follows from $(4.3.8)$ and $(4.3.9)$.
\par Next we consider the following equation
\par $(4.4.3)$ $(a_1,b_1)(a,b)=(e,e)$, where $a\in A$, $b\in B$.
\par From $(2.1.2)$ for fan loops $A$ and $B$, $(4.3.8)$ and $(4.3.9)$ we deduce that \par $(4.4.4)$
$a_1=e/a$,
\\ consequently, $b_1b^{(e/a)}\xi ((e/a,b_1),(a,b))=e$
and hence \par $(4.4.5)$ $b_1=e/[b^{(e/a)}\xi
((e/a,b^{(e/a)}),(a,b))]$. \\ Thus $a_1\in A$ and $b_1\in B$ given
by $(4.4.4)$ and $(4.4.5)$ provide a unique solution of $(4.4.3)$.
\par Similarly from the following equation
\par $(4.4.6)$ $(a,b)(a_2,b_2)=(e,e)$, where $a\in A$, $b\in B$
we infer that
\par $(4.4.7)$ $a_2=a\setminus e$, \\ consequently,
$bb_2^a\xi ((a,b),(a\setminus e,b_2))=e$ and hence \par $b_2^a=
b\setminus [\xi ((a,b),(a\setminus e,b_2))]^{-1}$ \\ by Conditions
$(2.1.1)$, $(2.1.2)$ and $(4.3.3)$ for fan loops $A$ and $B$. On the
other hand, $(b_2^a)^{e/a}=b_2\eta (e/a,a,b_2)$, consequently, by
Lemmas 2.2, 2.3 and the conditions of Section 4.3
\par $(4.4.8)$ $b_2= (b\setminus [\xi ((a,b),(a\setminus
e,(b\setminus e)^{e/a}))]^{-1})^{e/a})/\eta (e/a,a,(b\setminus
e)^{e/a}) $.
\\ Thus Formulas $(4.4.7)$ and $(4.4.8)$ provide a unique solution of
$(4.4.6)$.
\par Next we put $(a_1,b_1)=(e,e)/(a,b)$ and
$(a_2,b_2)=(a,b)\setminus (e,e)$ and
\par $(4.4.9)$ $(a,b)\setminus (c,d)=((a,b)\setminus (e,e))(c,d)
p((a,b),(a,b)\setminus (e,e), (c,d))$;
\par $(4.4.10)$ $(c,d)/(a,b)=[t((c,d),(e,e)/(a,b),(a,b))]^{-1}(c,d)((e,e)/(a,b))$
\par and $e_G=(e,e)$,
where $G=A\bigotimes ^{\phi , \eta , \kappa , \xi }B$. Therefore
Properties $(2.1.1)$-$(2.1.3)$ and $(2.1.9)$ are fulfilled for
$A\bigotimes ^{\phi , \eta , \kappa , \xi }B$.

\par {\bf 4.5. Definition.} The fan loop
$A\bigotimes ^{\phi , \eta , \kappa , \xi }B$ provided by Theorem
4.4 we call a smashed product of fan loops $A$ and $B$ with smashing
factors $\phi $, $\eta $, $\kappa $ and $\xi $.

\par {\bf 4.6. Corollary.} {\it Suppose that the conditions of Remark 4.3
are fulfilled and $A$ and $B$ are topological $T_1$ fan loops and
smashing factors $\phi $, $\eta $, $\kappa $, $\xi $ are jointly
continuous by their variables. Suppose also that $A\bigotimes ^{\phi
, \eta , \kappa , \xi }B$ is supplied with a topology induced from
the Tychonoff product topology on $A\times B$. Then $A\bigotimes
^{\phi , \eta , \kappa , \xi }B$ is a topological $T_1$ fan loop.}
\par {\bf 4.7. Corollary.} {\it If the conditions of Corollary 4.6
are satisfied and loops $A$ and $B$ are locally compact, then
$A\bigotimes ^{\phi , \eta , \kappa , \xi }B$ is locally compact.
Moreover, if $A$ and $B$ satisfy Condition $(3.5.2)$ and ranges of
$\eta $, $\kappa $, $\xi $ are contained in $N_0(A)N_0(B)$, then
$A\bigotimes ^{\phi , \eta , \kappa , \xi }B$ satisfies Condition
$(3.5.2)$.}
\par {\bf Proof.} Corollaries 4.6 and 4.7 follow immediately from
Theorems 2.3.11, 3.2.4, 3.3.13 in \cite{eng} and Theorem 4.4.

\par {\bf 4.8. Remark.} From Theorems 4.1, 4.4 and Corollaries
4.2, 4.6, 4.7 it follows that taking nontrivial $\phi $, $\eta $,
$\kappa $ and $\xi $ and starting even from groups with nontrivial
$N(G_j)$ or $N(A)$ and $G_j/\cdot/N(G_j)$ or $A/\cdot /N(A)$ it is
possible to construct new fan loops with nontrivial $N_0(G)$ and
ranges $t_G(G,G,G)$ and $p_G(G,G,G)$ of $t_G$ and $p_G$ may be
infinite and nondiscrete. With suitable smashing factors $\phi $,
$\eta $, $\kappa $ and $\xi $ and with nontrivial fan loops or
groups $A$ and $B$ it is easy to get examples of fan loops in which
$e/a\ne a\setminus e$ for an infinite family of elements $a$ in
$A\bigotimes ^{\phi , \eta , \kappa , \xi }B$.
\par {\bf 4.9. Conclusion.} The results of this article
can be used for further studies of measures on homogeneous spaces
and noncommutative manifolds related with loops. Besides
applications of left invariant measures on loops outlined in the
introduction it is interesting to mention possible applications in
mathematical coding theory and its technical applications
\cite{blautrctb,petbagsychrtj,srwseabm14}, because frequently codes
are based on topological-algebraic binary systems and measures.
Another very important applications are in representation theory of
loops and harmonic analysis on loops, mathematical physics, quantum
field theory, quantum gravity, gauge theory, etc.

\end{document}